\documentclass[11pt,a4paper]{article}
\usepackage[utf8]{inputenc} % This allows to type UTF-8 characters like Ä, Ü
\usepackage[T1]{fontenc}  % Tries to use Postscript Type 1 Fonts for better rendering
\usepackage{lmodern} % Provides the Latin Modern Font which offers more glyphs than the default Computer Modern
\usepackage{pdfpages}
\usepackage{scrhack}
\usepackage[main = english, ngerman]{babel} 
\usepackage{csquotes}
\usepackage{amsmath,amssymb,amsthm,amsfonts,bm}  % Provides all mathematical commands
\usepackage{mathtools}
\numberwithin{equation}{section}
\usepackage{tabularx}
\usepackage{booktabs} % Nice rules for tables. Usage \begin{tabular}\toprule ... \midrule ... \bottomrule
\usepackage{enumitem}
\usepackage{siunitx}  
\usepackage{setspace}
\usepackage{graphicx}
\usepackage{wrapfig}
\usepackage{subcaption}
\usepackage{hyperref}  % Provides clickable links in the PDF-document for \ref
\usepackage{bookmark}
\usepackage[noabbrev, capitalize]{cleveref}

\usepackage{xcolor}
\usepackage{grffile}
\usepackage{float} 
\usepackage{tikz, pgf}
\usetikzlibrary{fit,decorations.pathreplacing,patterns,shapes}
\usetikzlibrary{shapes,positioning,backgrounds}
\usetikzlibrary{plotmarks}
\usetikzlibrary{arrows.meta, automata}
\usetikzlibrary{calc,intersections}
\usetikzlibrary{external}
\usepackage{multirow}
\usepackage{pgfplots}
\pgfplotsset{compat=newest}
\usepgfplotslibrary{patchplots}

\usepackage{url}

\usepackage{a4wide} 
\usepackage{adjustbox}
\usepackage{titlesec}

\RequirePackage{ifthen}
\RequirePackage{listings} 

\usepackage{verbatim}

\theoremstyle{plain}
\newtheorem{theorem}{Theorem}[section] 
\newtheorem*{theorem*}{Theorem}
 
\newtheorem*{lemma*}{Lemma} 
 
\newtheorem*{korollar*}{Corollary} 
 
\newtheorem*{satz*}{Satz} 
\theoremstyle{definition}

\newtheorem*{definition*}{Definition} 

\newtheorem*{beispiel*}{Beispiel} 
\theoremstyle{remark} 
\newtheorem{remark}[theorem]{Remark}
\newtheorem*{remark*}{Remark}

\newtheorem*{anmerkung*}{Anmerkung}

\newcommand{\x}{\bm{x}}
\newcommand{\vv}{\bm{\mathrm{v}}}
\newcommand{\w}{\bm{\mathrm{w}}}
\newcommand{\outernormal}{\vec{\eta}}
\newcommand{\R}{\mathbb{R}}
\newcommand{\N}{\mathbb{N}}

\newcommand{\E}{\mathcal{E}}
\newcommand{\J}{\mathcal{J}}
\newcommand{\JJ}{\bm{\J}}
\newcommand{\dd}{\, \mathrm{d}}

\newcommand{\out}{{\mathrm{out}}}
\newcommand{\inn}{{\mathrm{in}}}

\newcommand{\numeric}{{\mathrm{num}}}
\newcommand{\diff}{{\mathrm{diff}}}
\newcommand{\adv}{{\mathrm{adv}}}

\renewcommand{\S}{{\mathcal{S}}}
\newcommand{\queue}{{\mathrm{queue}}}
\newcommand{\flow}{{\mathrm{flow}}}

\newcommand{\lwr}{{\mathrm{road}}}
\newcommand{\poll}{{\mathrm{poll}}}

\newcommand\logeq{\mathrel{\vcentcolon\Leftrightarrow}}

\newcommand{\p}{{\tilde{p}}}

\newcommand{\V}{{\bm{V}}}

\DeclareMathOperator*{\argmax}{arg\,max}

\definecolor{rwth}   {RGB}{  0  84 159}
\definecolor{rwth-75}{RGB}{ 64 127 183}
\definecolor{rwth-50}{RGB}{142 186 229}
\definecolor{rwth-25}{RGB}{199 221 242}
\definecolor{rwth-10}{RGB}{232 241 250}

\definecolor{rot}   {RGB}{204   7  30}
\definecolor{orange}   {RGB}{246 168   0}
\definecolor{lila}   {RGB}{122 111 172}
\definecolor{bordeaux-50}{RGB}{205 139 135}

\usepackage{verbatim}

\usepackage[backend=biber,style=numeric,giveninits=true,hyperref=true,isbn=false,url=true,doi=true,maxnames=10]{biblatex}
\begin{document}

    \title{Speed limits in traffic emission models using multi-objective optimization}
   \author{Simone G\"ottlich\footnotemark[1], \; Michael Herty\footnotemark[2], \; Alena Ulke\footnotemark[1]}

  \footnotetext[1]{University of Mannheim, Department of Mathematics, 68131 Mannheim, Germany (goettlich@uni-mannheim.de, ulke@uni-mannheim.de)}

    \footnotetext[2]{RWTH Aachen University, Department of Mathematics, 52056 Aachen, Germany (herty@igpm.rwth-aachen.de)}

    \date{ \today }

\maketitle

%%%%%%%%%%%%%%%%%%%%%%%%%%%%%%%%%%%%%%%%%%%%
%
%
\begin{abstract}
	\noindent 
	Climate change compels a reduction of greenhouse gas emissions, yet vehicular traffic still contributes significantly to the emission of air pollutants. Hence, in this paper we focus on the optimization of traffic flow while simultaneously minimizing air pollution using speed limits as controllable parameters. We introduce a framework of traffic emission models to simulate the traffic dynamic as well as the production and spread of air pollutants.
	We formulate a multi-objective optimization problem for the optimization of multiple aspects of vehicular traffic.
	The results show that multi-objective optimization can be a valuable tool in traffic emission modeling as it allows to find optimal compromises between ecological and economic objectives. 

\end{abstract}


\section{Introduction}
The societal impact of vehicular traffic has increased significantly over the last years and so have its drawbacks affecting the public's health, the climate, and the environment~\cite{DAmato2008,Epstein2005}. %~\cite{UBA}
As vehicular traffic emits air pollutants, it is directly linked to phenomena like air pollution, climate change, and global warming~\cite{Ramanathan2009}.

In this paper, we focus on the speed limit adjustments because (a) studies suggest that reducing speed limits reduces emissions~\cite{Panis2006}, and (b) this countermeasure is easy and inexpensive to implement, while also having the added benefits of reducing traffic noise in urban areas, cf.~\cite{Brink2022}, and improving road safety, cf.~\cite{Aljanahi1998}.

The ubiquity of vehicular traffic has motivated a vast branch of research, namely the traffic flow modeling. The models developed range from micro- to macroscopic, and from first to second order, cf.~\cite{Coclite2005,Garavello2006b,Garavello2006,Lighthill1955,Richards1956,Treiber2013} for an overview. In particular, much attention has been paid to economic aspects of traffic such as travel time, ramp metering, or traffic flow and their optimization, e.g., by controlling speed limits or traffic lights, cf.~\cite{annunziata2007optimization,Goatin2015,Goettlich2015,Treiber2013}.
There is also a growing interest in modeling traffic emissions~\cite{Panis2006} and monitoring concentrations of air pollutants, e.g., by simulating their spread due to advection and diffusion~\cite{Stoeki2011}. 
Recently, the idea of combining these approaches has been introduced in~\cite{Alvarez2017,Balzotti2022,berrone2012coupling,Pasquier2023,zegeye2011variable}.
Their common idea is to couple a traffic model with a model describing the dispersion of air pollutants. The coupling is achieved by estimating the emission of air pollutants based on quantities provided by the traffic flow model, and then, using this estimate as a source for the dispersion model.
The authors of~\cite{Alvarez2018} use this modeling approach to formulate an optimization problem aimed at finding the best extension of a road network to minimize traffic pollution.
Furthermore, there is an increasing interest in jointly optimizing economic and environmental aspects of traffic. For example,~\cite{VazquezMendez2019} follows a cooperative approach. Here, the authors minimize traffic pollution while simultaneously minimizing travel times and maximizing the outflow of the network by controlling the network management, i.e., the incoming capacity of the roads and the drivers' preferences. 
The same problem is studied in~\cite{GarciaChan2022}, where the authors instead follow a non-cooperative approach by formulating the problem as a bi-level Stackelberg game.
The approach of joint optimization of economic and environmental aspects is also used in different applications such as wastewater management, cf.~\cite{Alvarez2010,Alvarez2008}, or the management of industrial plants~\cite{Alvarez2015}.

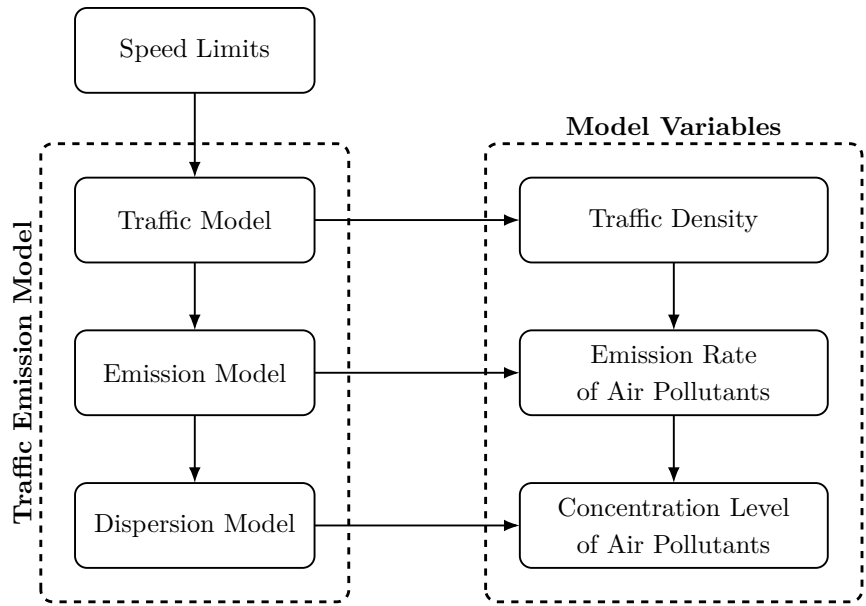
\begin{figure}[ht]
	\centering
	\scalebox{0.9}{
\begin{tikzpicture}[
	> = {Latex[scale=1.0]}
]
   	\pgfmathsetmacro{\hx}{3.5}
   	\pgfmathsetmacro{\hxRight}{4.5}
   	\pgfmathsetmacro{\hySpeed}{1.25}
   	\pgfmathsetmacro{\hy}{1.25}
   	\pgfmathsetmacro{\dx}{3}
   	\pgfmathsetmacro{\dy}{1}
   	\pgfmathsetmacro{\dh}{0.5}

	% box for speed limts on the very left
	\draw[rounded corners = 5pt, thick] 
   		(0, {3 *\hy + 2 * \dy + \hySpeed}) rectangle (\hx, {4 *\hy + 2 * \dy + \hySpeed}) 
   		node [midway] {Speed Limits};
   	\draw[->, thick] 
   			(\hx / 2, {3 *\hy + 2 * \dy + \hySpeed}) to (\hx / 2, {2 * (\hy + \dy) + \hy});

	% left column of boxes 
	\foreach \yy/\text in {%3/{Speed Limits},
						   2/{Traffic Model},
						   1/{Emission Model},
						   0/{Dispersion Model}}{
   		\draw[rounded corners = 5pt, thick] 
   			(0, {\yy * (\hy + \dy)}) rectangle (\hx, {\yy * (\hy + \dy) + \hy})
   				node [midway] {\text};
   	}
   	
   	% arrows pointing downwards (left)
   	\foreach \yy in {1,2}		
   		\draw[->, thick] 
   			(\hx / 2, {\yy * (\hy + \dy)}) to (\hx / 2, {\yy * (\hy + \dy) - \dy});
   	
   	% arrows pointing downwards (left)	
   	 \foreach \yy in {1,2}		
   		\draw[->, thick] 
   			(\hx + \hxRight / 2 + \dx, {\yy * (\hy + \dy)}) to (\hx + \hxRight / 2 + \dx, {\yy * (\hy + \dy) - \dy});

   	% right column of boxes
   	\draw[->, thick] (\hx, \hy/2) to (\hx + \dx,\hy/2);
   	\draw[rounded corners = 5pt, thick](\hx + \dx, 0) rectangle (\hx + \hxRight + \dx, \hy)
   		node [rectangle split, rectangle split parts = 2, midway] {Concentration Level \nodepart{second} of Air Pollutants};
   		
   	\draw[->, thick] 
   		(\hx, {(\hy + \dy) + \hy/2}) to (\hx + \dx, {(\hy + \dy) + \hy/2});
   	\draw[rounded corners = 5pt, thick] 
   		(\hx + \dx, {(\hy + \dy)}) rectangle (\hx + \hxRight + \dx, {(\hy + \dy) + \hy})
   		node [rectangle split, rectangle split parts = 2, midway] {Emission Rate \nodepart{second} of Air Pollutants};
   		
   	\draw[->, thick] 
   			(\hx, {2 * (\hy + \dy) + \hy/2}) to (\hx + \dx, {2 * (\hy + \dy) + \hy/2});
   	\draw[rounded corners = 5pt, thick] 
   		(\hx + \dx, {2 * (\hy + \dy)}) rectangle (\hx + \hxRight + \dx, {2 * (\hy + \dy) + \hy})
   		node [midway] {Traffic Density};

   	% surrounding boxes
   	\draw[rounded corners = 5pt, very thick, dashed]
   		(-\dh, - \dh) rectangle (\hx + \dh, {2 * \dy + 3* \hy + \dh});
 	%\node[rotate=0, above] at (\hx / 2, 3 * \hy + 2 * \dy + \dh) {\textbf{Simulation}};
 	\node[rotate=90, above] at (-\dh , {(3 / 2 * \hy + \dy}) {\textbf{Traffic Emission Model}};
   		
   	\draw[rounded corners = 5pt, very thick, dashed]
   		(\hx + \dx -\dh, - \dh) rectangle (\hx + \hxRight + \dx + \dh, {2 * \dy + 3* \hy + \dh});
   	 \node[rotate=0, above] at (\hx + \dx + \hxRight / 2, 3 * \hy + 2 * \dy + \dh) {\textbf{Model Variables}};
   	 %\node[rotate=90, above] at (\hx + \dx - \dh, {(3 / 2 * \hy + \dy}) {\textbf{Model Variables}};

\end{tikzpicture}
}
	\caption{Modeling framework of the traffic emission model including the dependence on the speed limits and
		the variables corresponding to the individual models.}
	\label{fig: model hierarchy}
\end{figure}

We extend this line of research by addressing the following question: Is it possible to maximize an economic aspect of traffic, such as traffic flow, while also minimizing the contribution to air pollution by controlling speed limits?
The novelty here lies in the approach to control the speed limits. We model this problem using multi-objective optimization.

%\begin{itemize}
First, we introduce the framework of speed-limit-dependent traffic emission models using a traffic, an emission, and a dispersion model where each depends on the quantities of the previous model. We illustrate this hierarchy in \cref{fig: model hierarchy}. The particular choice of these three models follows~\cite{Alvarez2017} with an additional explicit speed limit dependence.
Second, we formulate a multi-objective optimization problem aimed at maximizing traffic flow while simultaneously minimizing air pollution by controlling speed limits. 
Our approach to solve the proposed optimization uses the \emph{first-discretize-then-optimize} technique because the partial differential equations (PDEs) involved in the traffic emission model lack a closed-form solution in general.
%\end{itemize}

\medbreak
This article is structured as follows:  
First, we introduce the framework of so-called speed-limit-dependent traffic emission models which is split into the traffic model (\cref{sec: lwr}) and the emission and dispersion model (\cref{sec: emission model}).
Then, in \cref{sec: optimization}, we derive a multi-objective optimization problem capturing the goal of maximizing traffic flow while minimizing air pollution.
We cover the numerical treatment of our traffic emission model in \cref{sec: numerics} and conclude with a proof-of-concept example in \cref{sec: experiments}.

\section{Traffic flow dynamics}\label{sec: lwr}

The traffic dynamic on a road network consists of two components: the traffic dynamics on the roads, and the coupling at the junctions.
The traffic on the roads is described by a set of hyperbolic conservation laws, while the junctions are modeled by coupling conditions which encode boundary conditions at the end and beginning of each road~\cite{Coclite2005,Garavello2006b,Garavello2006}.

\medbreak
A road network is a directed graph \( (V, E) \) where each edge \( e \in E\) represents an unidirectional road associated with an interval \( I_e = [0,L_e^\lwr] \) and each vertex \(v \in V\) represents a junction.
A junction \(v\) can also be identified with a \((n+m)\)--tuple \( (e_1,\ldots, e_n,e_{n+1},\ldots,e_{n+m})\) where the first \(n\)--tuple corresponds to the incoming and the second \(m\)--tuple to the outgoing roads of \(v\), see also~\cite{Coclite2005,Garavello2006b,Garavello2006}.
We call a junction an internal junction if it has roads directed both into and away from the node. Conversely, we refer to a junction as an external junction, if it has either roads directed into or away from the node.
%In \cref{fig: sample network}, we provide an illustration for these definitions. 

%\begin{figure}[ht]%
%	\centering
%	\input{./tikz/square_network_small.tex}
%	\caption{An example of a road network. Internal junctions are indicated by (\CircleFilled) and e%xternal ones by (\Circle).}
%	\label{fig: sample network}
%\end{figure}

\subsection{The traffic dynamic on roads}
We describe the traffic dynamic on the roads of a network according to the Lighthill-Whitham-Richards (LWR) model~\cite{Lighthill1955, Richards1956}, which characterizes it in terms of the traffic density \( \rho_e(s,t)\) for given flux functions \( Q_e\) where \( s \in I_e\) and \(t \in \R^+\).
For each road \( e \in E \), the model reads
\begin{subequations}\label{eq: lwr}
	\begin{align}
		\partial_t \rho_e(s,t) + \partial_s Q_e(\rho_e(s,t), V_e^{\max}) & = 0, 
		&& \text{for } (s,t) \in I_e \times (0,T), \\
		\rho_e(s,0) & = \rho_e^0(s), && \text{for } s \in I_e.
	\end{align}
\end{subequations}
The flux functions \( Q_e\) are assumed to be dependent on the speed limits $V_e^{\max}$. The speed limits are constant over time and space, but may differ for every road.
We set the flux function according to Greenshields' model \cite{Greenshields1935}:
\begin{equation*}%\label{eq: Greenshields flux}
	Q_e(\rho, V_e^{\max}) \coloneqq V_e^{\max} \rho \left( 1 - \frac{\rho}{\rho_e^{\max}} \right),
\end{equation*}
where \( V_e^{\max} > 0 \) denotes the speed limit. 
The parameter \( \rho_e^{\max} \) describes the maximal density. In the Greenshields model, the critical density \( \rho_e^c \coloneqq \argmax_{0 \le \rho \le \rho_e^{\max}} Q_e(\rho, V_e^{\max})\) is independent of the speed limit policy as it is given by \( \rho_e^c = \rho_e^{\max} / 2 \). 
The traffic densities are affected by the full speed limit policy due to the interaction of traffic at junctions. We denote the speed limit policy on the network by the vector \( \V^{\max} = (V_e^{\max})_{e \in E}\) and hence, may write  \( \rho_e(s,t) = \rho_e(s,t,\V^{\max}) \).

Notice, that the Greenshields' model is only one possible option. The following considerations can be carried out for strictly concave flux functions.

\medbreak
The model~\labelcref{eq: lwr} is well-posed in the sense that it admits a (weak) solution for every road~\(e\) if additional boundary conditions are given. The claim follows directly from \cite{Garavello2006b,Garavello2006} as for fixed speed limits \( V_e^{\max} \), the model~\labelcref{eq: lwr} coincides with the \glqq classic\grqq{} LWR model.

\subsection{The traffic dynamic at junctions}\label{sec: junctions}
To characterize the traffic at junctions, we impose suitable boundary conditions.
In particular, we distinguish between two cases:
boundary conditions stemming from an external junction, where we prescribe in- or outflow rates, and boundary conditions stemming from an internal junction, where the flow rates arise by coupling conditions, see e.g.,~\cite{Coclite2005,Garavello2006b,Garavello2006,Goatin2015}.

%\medbreak 
At internal junctions, we suppose that the traffic dynamic obeys a physical condition, the conservation of mass. We follow~\cite{Coclite2005,Garavello2006b,Garavello2006} and consider additional rules at a junction besides the conservation of mass:
%\todo[noline]{\cite{Coclite2005,Garavello2006} doppelt zitiert. Einmal weg lassen?  Und Referenz zu \cite{Goatin2015} notwendig?}
\begin{enumerate}[label=(\roman*), ref=(\roman*)]
	\item At an internal junction, the traffic distributes from the incoming to the outgoing roads
	according to fixed, time-independent parameters that encode the drivers' preferences.\label{item: rule i}
	\item At a merging junction, where the number of incoming roads is greater than the number of 
	outgoing roads, fixed, time-independent parameters determine the priority among the incoming roads.\label{item: rule ii}
	\item The drivers aim to maximize the flux across the junctions with respect to 
	rule~\labelcref{item: rule i} or~\labelcref{item: rule ii}, respectively.
\end{enumerate}
These coupling conditions take the form of constrained maximization problems, one for each internal junction. Their solutions provide the in- and outflow at the boundary of the roads.

%We denote the flow arriving at the junction from the incoming road \( e_i\) by \( Q_{e_i}^\inn\) and the flow exiting the junction via the outgoing road \(e_j\) by \( Q_{e_j}^\out\).
%
To express the coupling conditions explicitly, we define the demand and supply function for each road \( e\):
\begin{equation*}
	D_e(\rho,V_e^{\max}) = \begin{cases} Q_e(\rho,V_e^{\max}), & \rho \le \rho_e^c, \\ Q_e^{\max}(V_e^{\max}), & \rho > \rho_e^c, \end{cases}
	\quad \text{and} \quad
	S_e(\rho,V_e^{\max}) = \begin{cases} Q_e^{\max}(V_e^{\max}), & \rho \le \rho_e^c, \\ Q_e(\rho,V_e^{\max}), & \rho > \rho_e^c. \end{cases}
\end{equation*}
and introduce the abbreviations: %[noline]{Ist sinnvolle bei den Abk\"urzungen \(\inn\) und \( \out \) zu tauschen?}
\begin{align*}
	D_i^\inn(t,\V^{\max}) &  \coloneqq D_{e_i}(\rho_{e_i}(L_{e_i}^\lwr,t,\V^{\max}),V_{e_i}^{\max}) && \text{for } i= 1, \ldots,n, \\
	S_j^\out(t,\V^{\max}) & \coloneqq S_{e_j}(\rho_{e_j}(0, t, \V^{\max}),V_{e_j}^{\max})
	&& \text{for } j = n+1, \ldots, n+m.
\end{align*}
Further, we introduce the following notation: The flow arriving at a junction from the incoming road \( e_i\) is denoted by \( Q_{e_i}^\inn\). Similarly, we denote flow exiting the junction via the outgoing road \( e_j\) by \( Q_{e_j}^\out\).
Then, the coupling condition at a junction at a fixed time \( t \in [0,T]\) takes the form: 
\begin{subequations}\label{eq: maximization at junction}
\begin{align}
	\max_{Q_{e_i}^\inn(t,\V^{\max})} & \,\,\, \sum_{i=1}^n Q_{e_i}^\inn(t,\V^{\max}) & \\
	\text{subject to} & \,\,\, \underbrace{\sum_{i=1}^n \alpha_{j,i} Q_{e_i}^\inn(t, \V^{\max})}_{\eqqcolon Q_{e_j}^\out(t,\V^{\max})}  \le S_j^\out(t,\V^{\max}), \\
	& \,\,\,  0 \le Q_{e_i}^\inn(t,\V^{\max}) \le D_i^\inn(t,\V^{\max}),
\end{align}
\end{subequations}
which is a linear optimization problem in finite dimensions.
The parameters \( \alpha_{j,i} \in (0,1) \) encode the drivers' preferences among the incoming roads and satisfy:
\(\sum_{j=n+1}^{n+m}\alpha_{j,i} = 1\).

The maximization problem~\labelcref{eq: maximization at junction} admits a solution for fixed speed limits given some additional assumptions, see e.g., \cite{Garavello2006b}, making the distribution of traffic at a junction a well-posed problem.
For an in-depth and rigorous discussion on coupling conditions at junctions, we refer to \cite{Garavello2006b} and the references therein.

\medbreak
In the following, we provide the explicit solution to the maximization problem~\labelcref{eq: maximization at junction} for the case of the one-to-one, the diverging one-to-two, and the merging two-to-one junction since our proof-of-concept example in \cref{sec: experiments} is solely based on these three junctions. The solutions are taken from \cite{Goatin2015}. 

We start with the one-to-one junction.
In this case, the flux across the junction is given the minimum of the supply of the outgoing and the demand of the incoming road:
\begin{equation}\label{eq: 1to1 junction}
	Q^\inn_{e_1}(t,\V^{\max}) = Q^\out_{e_2} (t,\V^{\max})
	= \min\{ D_1^\inn(t,\V^{\max}), S_1^\out(t,\V^{\max})) \}.
\end{equation}
%This conditions translates to: the incoming road sends at most its current demand but not more than the outgoing road can accommodate. 

For the diverging one-to-two junction the parameter \( \alpha_{j,i} \) describes the percentage of drivers arriving from road \( e_i\) at the junction, that take the outgoing road \(e_j \). The fluxes across such a junction is given by
\begin{subequations}\label{eq: 1to2 junction}
	\begin{align}
		Q^\out_{e_2}(t,\V^{\max}) & = \min\{ \alpha_{2,1} D_1^\inn(t,\V^{\max}), S_2^\out(t,\V^{\max}) \}, \\
		Q^\out_{e_3}(t,\V^{\max}) & = \min\{ \alpha_{3,1} D_1^\inn(t,\V^{\max}), S_3^\out(t,\V^{\max}) \}, \\
		Q^\inn_{e_1}(t,\V^{\max}) & = Q^\out_{e_2}(t,\V^{\max}) + Q^\out_{e_3}(t,\V^{\max}).
	\end{align}
\end{subequations}
The outflow of the incoming road is given by the sum of the inflow of the two outgoing roads, i.e., the number of vehicles entering the outgoing roads is the same as leaving the incoming road. Due to~\labelcref{eq: 1to2 junction} the number of vehicles is preserved.
%For the inflow the \( j \)-th outgoing road, we assume that \( \alpha_{j,1} \) percent of the drivers arriving at the junction, enter it. Thus, the incoming road sends \( \alpha_{j,1} \) percent of its demand to the outgoing road \( e_j^\out\).

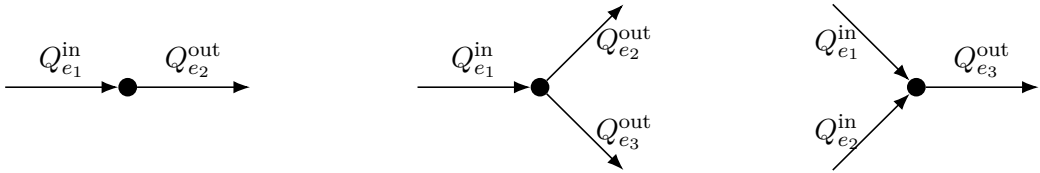
\begin{figure}[ht]
	\centering
	\begin{subfigure}[t]{.32\textwidth}
		\centering
		\scalebox{1}{
\begin{tikzpicture}[
	> = {Latex[scale=1.1]}, % arrow head style
    %shorten > = 1pt, % don't touch arrow head to node
    % auto,
    node distance = 1.75cm, % distance between nodes
    semithick % line style
]

	\tikzstyle{inner state}=[
    	fill=black, 
        circle,
        inner sep=0pt,
        outer sep=0pt,
        minimum size=7pt
   	]
        
     \tikzstyle{empty state}=[]

     \node[inner state] (1) {};        
     \node[empty state] (2) [right of=1] {}; 
     \node[empty state] (3) [left of= 1] {};      
     \node[empty state] (4) [below left of=1] {}; 
     \node[empty state] (5) [above left of= 1] {}; 

     \path[->] (1) edge node[midway, above] {$ Q^\out_{e_2} $} (2); %{1} (2);
     \path[->] (3) edge node[midway, above] {$ Q^\inn_{e_1} $} (1); %{3} (3);
\end{tikzpicture}
}
		\subcaption{One-to-one junction.}
	\end{subfigure} \hfill
	\begin{subfigure}[t]{.32\textwidth}
		\centering
		\scalebox{1}{
\begin{tikzpicture}[
	> = {Latex[scale=1.1]}, % arrow head style
    %shorten > = 1pt, % don't touch arrow head to node
    % auto,
    node distance = 1.75cm, % distance between nodes
    semithick % line style
]

	\tikzstyle{inner state}=[
    	fill=black, 
        circle,
        inner sep=0pt,
        outer sep=0pt,
        minimum size=7pt
   	]
        
     \tikzstyle{empty state}=[]

     \node[inner state] (1) {};        
     \node[empty state] (2) [left of=1] {}; 
     \node[empty state] (3) [above right of=1] {};  
     \node[empty state] (4) [below right of=1] {};      

     \path[->] (2) edge node[midway, above] {$ Q^\inn_{e_1}   $} (1); %{1} (2);
     \path[->] (1) edge node[midway, right] {$ Q^\out_{e_2}   $} (3); %{3} (3);
     \path[->] (1) edge node[midway, right] {$ Q^\out_{e_3}   $} (4); %{3} (3);
\end{tikzpicture}
}
		\subcaption{One-to-two junction.}
	\end{subfigure} \hfill
	\begin{subfigure}[t]{.32\textwidth}
		\centering
		\scalebox{1}{
\begin{tikzpicture}[
	> = {Latex[scale=1.1]}, % arrow head style
    %shorten > = 1pt, % don't touch arrow head to node
    % auto,
    node distance = 1.75cm, % distance between nodes
    semithick % line style
]

	\tikzstyle{inner state}=[
    	fill=black, 
        circle,
        inner sep=0pt,
        outer sep=0pt,
        minimum size=7pt
   	]
        
     \tikzstyle{empty state}=[]

     \node[inner state] (1) {};        
     \node[empty state] (2) [right of=1] {}; 
     \node[empty state] (3) [above left of=1] {};  
     \node[empty state] (4) [below left of=1] {};         

     \path[->] (1) edge node[midway, above] {$Q^\out_{e_3}$} (2); %{1} (2);
     \path[->] (3) edge node[midway, left] {$Q^\inn_{e_1} $} (1); %{3} (3);
     \path[->] (4) edge node[midway, left] {$Q^\inn_{e_2}$} (1); %{3} (3);
\end{tikzpicture}
}
		\subcaption{Two-to-one junction.}
	\end{subfigure}
	\caption{In- and out flow rates at the three basis junctions.}
	\label{fig: basic junctions}
\end{figure}

At a merging two-to-one junction the parameters \( \beta_{j,i}\) regulate the priority of the incoming roads.
The parameter \( \beta_{3,i} \) describes the percentage of drives coming from the incoming road \( e_i \) that are allowed to take the only outgoing road:
\begin{subequations}\label{eq: 2to1 junction}
	\begin{align}
		Q_{e_1}^\inn(t,\V^{\max}) & = \min\{ D_1^\inn(t,\V^{\max}), 
			\gamma_1 \}, \\
		\text{where} \quad  \gamma_1 & \coloneqq  \max\{\beta_{3,1} S_1^\out(t,\V^{\max}), S_1^\out(t,\V^{\max}) - D_2^\inn(t,\V^{\max})\}, \\
		Q_{e_2}^\inn(t,\V^{\max}) & = \min\{ D_2^\inn(t,\V^{\max}), 
		 \gamma_2 \}, \\
		\text{where} \quad  \gamma_2 & \coloneqq \max\{\beta_{3,2} S_1^\out(t,\V^{\max}), S_1^\out(t,\V^{\max}) - D_1^\inn(t,\V^{\max}) \}, \\
		Q_{e_3}^\out(t,\V^{\max}) & = Q_{e_1}^\inn(t,\V^{\max}) 
		+ Q_{e_2}^\inn(t,\V^{\max}). 
	\end{align}
\end{subequations}

\medbreak
Next, we discuss the boundary conditions at external junctions. We distinguish between two cases: roads where vehicles exit the network, and roads where vehicles enter the network.
In the first case, we assume free flow and in the second, we follow \cite{Goatin2015,Herty2009} and prescribe an inflow rate \(q_e^\inn\) and model a queue in case the desired inflow exceeds available capacity of the road. 
The queue can be interpreted as an infinitely long road with no spatial extension, and allows only feasible flow to be sent onto the network via the access roads.
The modeling of a queue can redistribute the inflow over time, allowing the flow to enter the network at a later point in time at which the capacity allows it. Further, we are interested in sending the maximal possible inflow.
Thus, the evolution of the queue length \( \ell_e \) obey the following ordinary differential equation (ODE):
\begin{equation}\label{eq: queue ODE}
	\dot{\ell}_e(t) = q_e^\inn(t)- Q_e^\out(t,\V^{\max}), \quad
	\ell_e(0) = \ell_e^0,
\end{equation}
The flow rate \( Q_e^\out \) describes the outflow from the queue onto the road and is therefore speed-limit-dependent. It is determined by the minimum between the capacity of the road, to ensure feasible flow rates, and the demand \( D_e^\queue\) of the queue:
\begin{equation}\label{eq: queue inflow}
	Q_e^\out(t,\V^{\max}) = \min \bigg \{ D_e^\queue(t), S_e(\rho_e(0,t,\V^{\max}), V_e^{\max})\bigg  \}.
\end{equation}
The demand of the queue is composed of the prescribed inflow rate \(q_e^\inn\) and the flow rate \(q_e^\queue\) that describes the flow required to clear the current queue connected to road \(e\) immediately, i.e.,
\begin{equation}\label{eq: queue demand}
	D_e^\queue(t) = q_e^\inn(t) + q_e^\queue(t).
\end{equation}
The existence of a solution to the ODE~\labelcref{eq: queue ODE} follows from \cite{Garavello2006b} and the references therein. %Herty2009}.
Again, we want to highlight that the evolution of the queue length depends on the fluxes on the roads, making it speed-limit-dependent and we may write \( \ell_e(t) = \ell_e(t, \V^{\max})\).

\section{Modeling air pollution caused by vehicular traffic}\label{sec: emission model}
In this section, we cover the emission and dispersion model which simulate the traffic's contribution to air pollution. 
%As the contribution to air pollution is far too complex to depict in a mathematical model rigorously, we simplify the problem by taking the concentration level of air pollutants and their absolute amount in a control area as an indicator.
%These air pollutants include carbon mon- and dioxide, nitrogen oxides, or particulate matter.
%
%We approach the ansatz in two steps:
Our modeling approach is the following:
\begin{enumerate}[label=(\roman*), ref=(\roman*)]
	\item The emission model estimates the emission rate of air pollutants on a given
		road based on the quantities the traffic model provides.
		\label{item: emission model}
	\item The dispersion model describes the spread of the pollutants in a 
		two-dimensional domain \( \Omega \) representing the simulated city.
		We include traffic's emission of pollutants by employing the estimate
		from~\labelcref{item: emission model}.%   , cf. \cite{Alvarez2018,SecondOrderModel}.
\end{enumerate}

%We stress that the emission rate and the concentration level of the air pollutant are speed-limit-dependent because the flux functions of the traffic model explicitly depend on the speed limits.
%Thus, the traffic densities inherit this property which makes the emission rate speed-limit-dependent, too, as it is estimated based on the quantities the traffic model provides.
%Further, the spread of the pollutant depends on the speed limits since it does depend on the emission rate.
		
\subsection{The emission model}
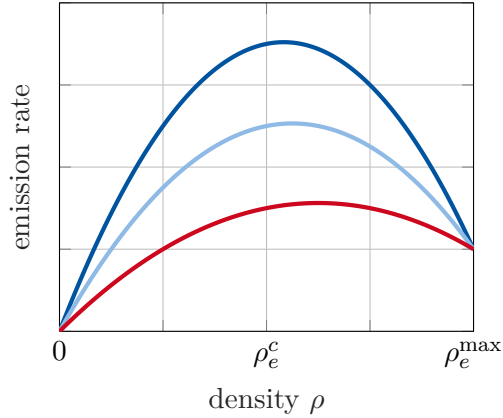
\begin{figure}[ht]
	\centering
	\scalebox{1.0}{
\begin{tikzpicture}
\begin{axis}[
width=0.35\textwidth,
height=0.28\textwidth,
at={(0\textwidth,0\textwidth)},
scale only axis,
xmin = 0, 
xmax = 1,
xtick={0, 0.25, 0.5, 0.75, 1},
xticklabels={0, , $\rho_e^c$, , $\rho_e^{\max}$},
xlabel style={font=\color{white!15!black}},
xlabel={density $\rho$},
ymin = 0, 
ymax = 0.5,
ytick={0, 0.125, 0.25, 0.375, 0.5},
yticklabels={, , , ,},
yticklabels={, , , ,},
ylabel style={font=\color{white!15!black}},
ylabel={emission rate},
%ticklabel style = {font=\normalsize},
axis background/.style={fill=white},
xmajorgrids,
ymajorgrids]
    \addplot[
        domain = 0:1,
        samples = 200,
        smooth,
        ultra thick,
        %densely dashdotted,
        rwth,
    ] {1.5 * x * ( 1 - x ) + x / 8};
    
    \addplot[
        domain = 0:1,
        samples = 200,
        smooth,
        ultra thick,
        rwth-50,
    ] {x * ( 1 - x ) + x / 8};
    
    \addplot[
        domain = 0:1,
        samples = 200,
        smooth,
        ultra thick,
        %densely dotted,
        rot,
    ] {0.5 * x * ( 1 - x ) + x / 8};

\end{axis}
\end{tikzpicture}
}
	\caption{Influence of the speed limit on the emission rate~\labelcref{eq: emission rate}
		for a fixed parameter \( \theta > 0 \). The corresponding speed limits are ordered as
		follows: \(\textcolor{rot}{{} V_e^{\max}} < \textcolor{rwth-50}{{} V_e^{\max}}	< \textcolor{rwth}{{} V_e^{\max}}\)}
	\label{fig: emission rate}
\end{figure}

For the estimation of the emission rate, we follow \cite{Alvarez2018} and use a linear combination of the traffic density and traffic flow. 
Therefore, we introduce a weighting parameter \( \theta \ge 0 \) which represents the relative influence of the traffic density to the emission rate compared to the traffic flow. 
%This parameter is chosen in dependence on the air pollutants under investigation.
Also considering the speed-limit dependence of the traffic density and flux function, the estimate of the emission rate \(\xi_e\) on road \( e \) reads:
\begin{equation}\label{eq: emission rate}
	\xi_e(s,t,\V^{\max}) = Q_e(\rho_e(s,t,\V^{\max}),\V^{\max}) + \theta \rho_e(s,t,\V^{\max})
		\quad \text{for } (s,t) \in I_e \times (0,T).
\end{equation}
For a fixed weighting parameter \( \theta > 0 \), we visualize the emission rate of road \( e \) in dependence on the traffic density in \cref{fig: emission rate}.

\begin{figure}[ht]
	\centering
	\begin{subfigure}[t]{.32\textwidth}
		\centering
		\scalebox{0.9}{
\begin{tikzpicture}[
	> = {Latex[scale=1.1]},
    very thick,
]
        
\pgfmathsetmacro{\L}{5}

% road network
\tikzstyle{state}=[
    fill=rot, 
    circle,
    inner sep=0pt,
    outer sep=0pt,
    minimum size=8pt
]
        
\draw[black!30!white, fill=black!3!white, line width=1.25pt] (0, 0) rectangle (\L, \L); 
			
\path[->] (\L / 2,0) edge node[midway, left] {$\sigma_1(I_1)$} (\L / 2, \L / 2);
%\node[state] at (\L / 2, \L / 8) {};
%\node[above right] at (\L / 2, \L / 8) {$\textcolor{rot}{\sigma_1(s)}$};
\path[->] (\L / 2, \L / 2) edge node[midway, above] {$\sigma_2(I_2)$} (\L, \L / 2);
	
\node[below right] at (0, \L) {$\textcolor{black!30!white}{\bm{\Omega}}$};

\end{tikzpicture}
}
		\subcaption{Courses of the roads parameterized by the curves \( \sigma_1 \) and \( \sigma_2\).}
		\label{fig: multiple roads a}
	\end{subfigure} \hfill
	\begin{subfigure}[t]{.32\textwidth}
		\centering
		\scalebox{0.9}{
\begin{tikzpicture}[
	> = {Latex[scale=1.1]},
    very thick,
]
        
\pgfmathsetmacro{\L}{5}
%\pgfmathsetmacro{\h}{0.25}
%\pgfmathsetmacro{\r}{5.5}
\pgfmathsetmacro{\w}{1.25}

% road network
\tikzstyle{state}=[
    fill=rot, 
    circle,
    inner sep=0pt,
    outer sep=0pt,
    minimum size=8pt
]

\draw[black!30!white, fill=black!3!white, line width=1.25pt] (0, 0) rectangle (\L, \L); 

\draw[fill=rwth-10, draw=black!30!white] (\L / 2 - \w, 0) rectangle (\L / 2 + \w, \L / 2);
%\draw[fill=rwth-10, draw=black!30!white] (\L / 2, \L / 2 - \w) rectangle (\L, \L / 2 + \w);%
%\draw[fill=rwth-10, draw=rwth-10] (\L / 2, \L / 2 - \w) rectangle (\L / 2 + \w, \L / 2);
			
\path[->] (\L / 2, 0) edge node[midway, left] {} (\L / 2, \L / 2);
\path[->] (\L / 2, \L / 2) edge node[midway, above] {} (\L, \L / 2);

\node[state] at (\L / 2, \L / 8) {};
\node[above right] at (\L / 2, \L / 8) {$\textcolor{rot}{\sigma_1(s)}$};
\node[left] at (\L / 2 - \w, \L / 8) {$\textcolor{rot}{\mathcal{R}_1(s)}$};
\draw[rot] (\L / 2 - \w, \L / 8) -- (\L / 2 + \w, \L / 8);	

\node[below right] at (0, \L) {$\textcolor{black!30!white}{\bm{\Omega}}$};

%	\draw[fill=rwth-10, draw = none]
%		(\L - \r - \w,0) arc(180:90:{\r + \w}) -- (\L, \r + \w) 
%			-- (\L, \r - \w) arc(90:180:{\r - \w})
%			-- (\L - \r + \w,0)  -- (\L - \r - \w,0) -- cycle;
		 
%	\draw[-] (\L - \r - \w,0) arc(180:90:{\r + \w}) -- (\L, \r + \w);
%	\draw[-] (\L - \r + \w,0) arc(180:90:{\r - \w}) -- (\L, \r - \w);			
%	\draw[->] (\L - \r,0) arc(180:91:\r) -- (\L, \r);
	
%	\draw[-, rot] ({(\r + \w) * cos(157.5) + \L}, {(\r+\w) * sin(157.5)}) 
%		--  ({(\r - \w) * cos(157.5) + \L}, {(\r-\w) * sin(157.5)});
%	\node[state, fill=rot, draw=rot] at ({\r * cos(157.5) + \L}, {\r * sin(157.5)}) {};
	
%	\node[right] at ({\r * cos(157.5) + \L +0.25}, {\r * sin(157.5)}) 
%		{$\textcolor{rot}{\sigma_e(s)}$};
%	\node[right] at ({(\r - \w) * cos(157.5) + \L}, {(\r-\w) * sin(157.5)}) 
%		{$\textcolor{rot}{\mathcal{R}_e(s)}$};
%	
	
%	\node[below right] at (0, \L) {$\textcolor{black!30!white}{\bm{\Omega}}$};
%	\draw[black!30!white, fill=none, line width=1.25pt] (0, 0) rectangle (\L, \L);

\end{tikzpicture}
}
		\subcaption{Visualization of \( \mathcal{R}_e(s)\) for road \( e = 1\).}
		\label{fig: multiple roads b}
	\end{subfigure} \hfill
	\begin{subfigure}[t]{.32\textwidth}
		\centering
	    \scalebox{0.9}{
\begin{tikzpicture}[
	> = {Latex[scale=1.1]}, % arrow head style
    %shorten > = 1pt, % don't touch arrow head to node
    % auto,
    node distance = 4cm, % distance between nodes
    very thick % line style
]

\pgfmathsetmacro{\L}{5}
\pgfmathsetmacro{\w}{1.25}

\draw[black!30!white, fill=black!3!white, line width=1.25pt] (0, 0) rectangle (\L, \L); 

\draw[fill=rwth-10, draw=black!30!white] (\L / 2 - \w, 0) rectangle (\L / 2 + \w, \L / 2);
\draw[fill=rwth-10, draw=black!30!white] (\L / 2, \L / 2 - \w) rectangle (\L, \L / 2 + \w);
\draw[fill=rwth-50, draw=black!30!white] (\L / 2, \L / 2 - \w) rectangle (\L / 2 + \w, \L / 2);
			
\path[->] (\L / 2, 0) edge node[midway, left] {} (\L / 2, \L / 2);
\path[->] (\L / 2, \L / 2) edge node[midway, above] {} (\L, \L / 2);

\node[below right] at (0, \L) {$\textcolor{black!30!white}{\bm{\Omega}}$};

%	\draw[draw=black!30!white, fill=black!3!white] (0,-\L) rectangle (3*\L, 2*\L);
%	\node[below right] at (0, 2*\L) {$\textcolor{black!30!white}{\bm{\Omega}}$};
%			
%%	% roads with width > 0
	%\draw[fill=rwth-10, draw=black!30!white]  (0, -\w/2) rectangle (\L, \w/2);
	%\draw[fill=rwth-10, draw=black!30!white] 
%		(\L, - \w/2) rectangle (2*\L, \w/2);
%	\draw[fill=rwth-10, draw=black!30!white] 
%		(\L - \w/2, 0) rectangle (\L + \w / 2, \L);
%	\draw[fill=rwth-10, draw=black!30!white] 
%		(2*\L - \w/2, 0) rectangle (2 * \L + \w / 2, \L);
%	\draw[fill=rwth-10, draw=black!30!white] 
%		(\L, \L - \w/2) rectangle (2 * \L, \L + \w/2);
%	\draw[fill=rwth-10, draw=black!30!white] 
%		(2*\L, \L - \w/2) rectangle (3 * \L, \L + \w/2);
		
	% overlapping parts of the roads
%	\fill[rwth-50] (\L - \w / 2, 0) rectangle (\L + \w/2, \w/2);
%	\fill[rwth-50] (2*\L - \w / 2, 0) rectangle (2*\L, \w/2);
%	\fill[rwth-50] (2*\L - \w / 2, \L - \w/2) rectangle (2*\L + \w/2, \L);
%	\fill[rwth-50] (\L, \L - \w /2) rectangle (\L + \w/2, \L);

	% the network
%	\draw[->, black] (0,0) to (\L,0);
%	\draw[->, black] (\L,0) to (2*\L,0);
%	\draw[->, black] (\L,0) to (\L,\L);
%	\draw[->, black] (\L,\L) to (2*\L,\L);
%	\draw[->, black] (2*\L,\L) to (3*\L,\L);
%	\draw[->, black] (2*\L,0) to (2*\L,\L);
			
\end{tikzpicture}
}
		\subcaption{Courses of the roads including the width \( w_1 = w_2 \).}
		\label{fig: multiple roads c}
	\end{subfigure}
    \caption{Points belonging to the grey shaded area are not associated to any road, while
		points belonging to the light blue shaded area are associated to a single road, 
    	and points belonging to dark blue shaded area are associated to two roads simultaneously.}
	\label{fig: multiple roads}
\end{figure}
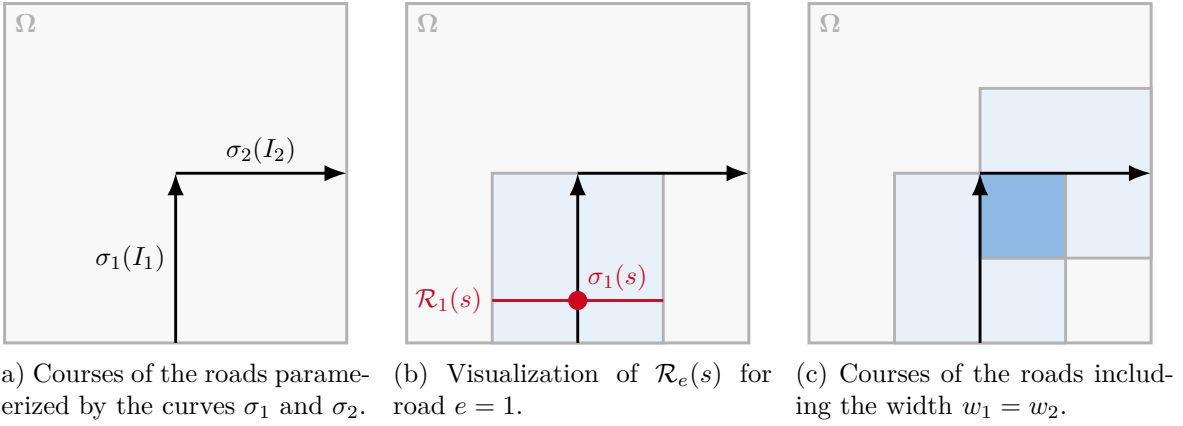

Furthermore, we assume that the road network represented by the graph \( (V,E)\), is contained in a domain \( \Omega \subset \R^2 \). To model the spread of air pollutants, which are emitted by vehicular traffic on the roads of a network, we introduce a source term that defines the emission rate at each point in \( \Omega \).  
Therefore, we parameterize each edge \(e \in E\) by a curve \( \sigma_e \colon I_e \rightarrow \Omega \).
%\begin{equation}\label{eq: curve}
%	\sigma_e \colon I_e \rightarrow \Omega, \quad s \mapsto \sigma_e(s) \in \Omega,
%\end{equation}
Its starting and end point are given by \( \sigma_e(0) \) and \( \sigma_e(L_e^\lwr) \), respectively to preserve the direction of \( e\). See \cref{fig: multiple roads a} for an illustration.

As the curves \( \sigma_e\) allow to identify a point \( s \in I_e \) with a point \( \x \in \Omega \), we map the values of the emission rates \( \xi_e \) to the two-dimensional domain \( \Omega \). 
Further, we assign a positive width \( w_e \) to each road and distribute the emission of air pollutants along the width \( w_e\).
Then the emission rate of road \( e \) within the domain \( \Omega \) is:
\begin{equation}\label{eq: emission rate in omega}
	\tilde{\xi}_e (\x,t,\V^{\max}) = \begin{cases}
		\frac{1}{w_e} \xi_e(s,t,\V^{\max}), & \text{if } \x \in \mathcal{R}_e(s), \\
		0, & \text{otherwise},
		\end{cases}
\end{equation}
where the set \( \mathcal{R}_e(s) \) describes the points along the width of the road \( e\):
\begin{equation*}
	\mathcal{R}_e(s) \coloneqq \bigg \{ \x \in \Omega \mathrel{\Big|}
		\Vert \sigma_e(s) - \x \Vert_2 \le \frac{w_e}{2} 
		\text{ and } \langle \sigma^\prime_e(s), \sigma_e(s) - \x \rangle = 0 \bigg \}.
\end{equation*}
In \cref{fig: multiple roads b}, we provide a visualization of this set for a road which is parallel to the \(y\)-axis.
This example reveals one problem: Based on our modeling approach, for two roads \( e \neq \tilde{e}\) the following relation can hold:
\begin{equation*}
	\mathcal{R}_e(s) \cap \mathcal{R}_{\tilde{e}}(\tilde{s}) \neq \emptyset
		\quad \text{where} \quad
		 s\in I_e \text{ and } \tilde{s} \in I_{\tilde{e}}.
\end{equation*}
This relation implies that there might exist points \( \x \in \Omega \) to which multiple emission rates are assigned, see \cref{fig: multiple roads c} for an illustration. 
This phenomenon also occurs when considering multiple, pairwise different roads, e.g., at junctions where more than two roads meet.

We resolve this problem by taking the average of the emission rates which are assigned to a point \( \x \) based on the model given in \cref{eq: emission rate in omega}.
%One can also think of other ways to resolve the ambiguity of the emission rate. However, this is subject to future research as this would exceed the scope of this work.
Therefore, we define the set \( \E(\x) \) that contains the roads to which the point \( \x \) can be associated given \( \{\sigma_e\}_{e \in E} \) and the width' \( \{w_e\}_{e \in E} \) of the roads:
\begin{equation*}
	 \E(\x) \coloneqq \{ e \in E \mid \x \in  \bm{\mathcal{R}}_e \}
	 	\quad \text{with} \quad 
	 \bm{\mathcal{R}}_e \coloneqq \bigcup_{s \in I_e} \mathcal{R}_e(s).
\end{equation*}
Thus, \( \vert \E(\x) \vert \) describes the number of roads assigned to \( \x \). 
Finally, the joint emission rate defined on \( \Omega \), which is called \emph{emission model} in the following, is given by:
\begin{equation}\label{eq: emission model}
	\xi(\x,t,\V^{\max}) \coloneqq \begin{cases}
		\frac{1}{\vert \E(\x) \vert } 
		\sum_{e \in \E(\x)} \tilde{\xi}_e(\x,t,\V^{\max}), & \text{if } \E(\x) \neq \emptyset, \\
		0, & \text{otherwise}.
	\end{cases}		
\end{equation}

\subsection{The dispersion model}\label{sec: dispersion model}
After having established an estimate for the emission of air pollutants, we now introduce the dispersion model to simulate the spread of these pollutants within the domain \( \Omega \) due to wind and horizontal diffusion. 
 follow the ideas of \cite{Alvarez2017,Balzotti2022,Parra2003,Skiba2000,Stoeki2011} and invoke an advection-diffusion equation with suitable boundary conditions and make two assumptions: vehicular traffic is the only source of pollution in our model, i.e., pollutants cannot enter the domain from outside, but they can leave it due to wind and diffusion.

To formally represent these assumptions, we introduce the in- and outflow boundary of \( \Omega \) for a given wind field. A wind field is a space- and time-dependent vector field \( \vv = (v_x, v_y)^\top\) and describes the wind dynamic in \( \Omega\). 
Additionally, we assume that \( \Omega \) is bounded, and its unit outward normal is denoted by \( \outernormal\).
%
%We suppose that the wind field is independent of time.
Then, the inflow boundary is defined by
\begin{equation*}
%	\Gamma^-(\vv) \times (0,T) 
%		\quad \text{with} \quad 
%		\Gamma^-(\vv) \coloneqq \{ \x \in \partial \Omega 
%			\mid \vv(\x) \cdot \outernormal(\x) < 0 \},
	S^-(\vv) \coloneqq \{ (\x,t) \in \partial \Omega \times [0,T] \mid \vv(\x,t) \cdot \outernormal(\x) < 0\},
\end{equation*}
which describes the boundary through which pollutants enter the domain from outside due to the wind dynamic.
Conversely, the outflow boundary describes the boundary where pollutants leave the domain due to the wind dynamic. Therefore, it is given by
\begin{equation*}
%	\Gamma^+(\vv) \times (0,T)
%		\quad \text{with} \quad 
%		\Gamma^+(\vv) \coloneqq \{ \x \in \partial \Omega 
%			\mid \vv(\x) \cdot \outernormal(\x) \ge 0 \}.
	S^+(\vv) \coloneqq \{ (\x,t) \in \partial \Omega \times [0,T] \mid \vv(\x,t) \cdot \outernormal(\x) \ge 0\}.
\end{equation*}	
Finally, given an initial concentration distribution \( \phi_0 \) of the air pollutants in \( \Omega \), the following advection-diffusion equation describes the evolution of the concentration \( \phi \) over time:
\begin{subequations}\label{eq: dispersion model}
	\begin{align}
		\partial_t \phi(\x,t) - \mu \Delta \phi(\x,t) 
			+ \vv \cdot \nabla \phi(\x,t) + \kappa \phi(\x,t)
			  & = \xi(\x,t,\V^{\max}),
			&& (\x,t) \in \Omega \times (0,T), \label{eq: advection diffusion}\\
		\mu \nabla \phi(\x,t) \cdot \outernormal(\x) - \phi(\x,t) \vv \cdot \outernormal(\x)  & = 0, &&  (\x,t) \in S^-(\vv), \label{eq: Dirichlet boundary} \\
			\nabla \phi(\x,t) \cdot \outernormal(\x) & = 0, 
		&&   (\x,t) \in S^+(\vv),  \label{eq: Neumann boundary} \\
		\phi(\x,0) & = \phi_0(\x), &&  \x \in \Omega, \label{eq: initial condition}
	\end{align}
\end{subequations}
where the constant \( \mu > 0 \) describes the diffusion coefficient of the air pollutants and the constant \( \kappa \) their extinction rate to account for (chemical) interactions with other substances in the air.
The boundary condition~\labelcref{eq: Neumann boundary}  allows air pollutants to stream over the outflow boundary, whereas the boundary condition~\labelcref{eq: Dirichlet boundary} guarantees that the combined diffusive and advective flow is zero on the inflow boundary. The later condition enforces that pollutants cannot enter from outside.

%\begin{remark} 
%	The speed limits enter the flux functions of the traffic model explicitly making the traffic density also dependent on the speed limits and likewise the emission rate \( \xi \).
%	Thus, the concentration \( \phi \) is also speed-limit-dependent because the source term of the dispersion model, the emission rate, depend on the speed limits.
%	Again, we can write \( \phi(\x,t) = \phi(\x,t,\V^{\max}) \).
%\end{remark}

The existence and uniqueness of a weak solution to the dispersion model~\labelcref{eq: dispersion model} follows from~\cite[Theorem 1]{Skiba2000}, which makes it a well-posed problem.

\section{The Multi-objective optimization approach}\label{sec: optimization}
After introducing the traffic emission model in \cref{sec: lwr,sec: emission model}, our next goal is to maximize the economic efficiency of vehicular traffic while also minimizing the environmental impact by adjusting the speed limits on the roads.
The formulation of the corresponding optimization problem is divided into two steps:
\begin{enumerate}[label=(\roman*), ref=(\roman*)]
	\item The description of \emph{economic efficiency} and the \emph{environmental impact} by
	scalar-valued and speed-limit-dependent functions. Here, we use ideas already established in the literature. \label{item: objectives}
	\item The derivation of a multi-objective optimization problem where the functions obtained from \labelcref{item: objectives} are optimized simultaneously. The approach of optimizing two aspects of traffic at the same time by controlling speed limits is unique in the traffic literature, at least to the best of our knowledge.
\end{enumerate}

\subsection{Quantification of economic and ecological aspects}\label{sec: objectives}
The main concern of this section is the discussion on indicators in the form of scalar-valued, speed-limit-dependent functions to measure the economic and ecological aspects of vehicular traffic. 
The idea is to utilize the quantities provided by the traffic emission model, which include the traffic density and flux, or the concentration level of air pollutants.

\medbreak
The economic efficiency of vehicular traffic can be quantified by various aspects -- for example, the travel time or the traffic flow, cf. \cite{Goettlich2015,Treiber2013}.
We maximize traffic flow by considering \emph{accumulated traffic flow}:
\begin{equation}\label{eq: J flow}
	\J_\flow(\V^{\max}) \coloneqq \sum_{e \in E} \int_0^T \int_{I_e} 
		Q_e(\rho_e(s,t, \V^{\max}), V_e^{\max}) \dd s \dd t.
\end{equation}

\medbreak
To determine the influence of vehicular traffic on the environment we measure its contribution to air pollution. Our modeling admits two sources:
The \emph{active} traffic on the road network and the \emph{idle} traffic in the queues.
The vehicles waiting in queues to enter the road network can contribute to air pollution since their engines may remain idle instead of being turned off.
We measure the air pollution caused by the active traffic on the network by the average mass of air pollutant in a given control area \( \Omega \) until a final time \( T > 0 \), cf. \cite{Alvarez2018}:
\begin{equation}\label{eq: J_diff original}
	\J_\diff(\V^{\max}) \coloneqq \frac{1}{T \vert \Omega \vert} 
		\int_0^T \int_{\Omega} \phi(\x,t, \V^{\max}) \dd \x \dd t,
\end{equation}
where \( \vert \Omega \vert \) denotes the area of \( \Omega \).
We measure the air pollution caused by idle traffic by the average amount of air pollutants emitted by idling vehicles waiting in queues. Assuming that all vehicles have the same emission rate when idle, we compute the pollution caused by idling vehicles based on the average queue length:
\begin{equation} 
	\J_\queue(\V^{\max}) \coloneqq \frac{1}{T} 
		\sum_{e \in E_\inn} \int_0^T \ell_e(t, \V^{\max}) \dd t,
\end{equation}
where set \( E_\inn \) contains all access roads of the road network, i.e., roads at which we impose an external inflow rate \( q_e^\inn\).
Thus, an indicator for the influence of vehicular traffic on the environment based on air pollution is:
\begin{equation} \label{eq: J eco}
	\J_\poll(\V^{\max}; \delta) \coloneqq \J_\diff(\V^{\max}) + \delta  \J_\queue(\V^{\max}),
\end{equation}
where \( \delta \ge 0\) describes the emission of air pollutants per time unit of an idle engine.

Besides the expression~\labelcref{eq: J_diff original} for \( \J_\diff\), we can employ another, equivalent one that allows us to solve a suitable PDE once, and then requires only a re-evaluation of the traffic model for each new speed limit policy, eliminating the need to solve the dispersion model.
The derivation of the alternative expression is inspired by \cite{Alvarez2018} and involves integration by parts and some adjoint calculus. The result reads:
\begin{equation}\label{eq: alternative expression J_diff}
	\J_\diff(\V^{\max}) = \int_0^T \int_{\Omega} \xi(\x,t,\V^{\max}) p(\x, t) \dd \x \dd t
	+ \int_{\Omega} \phi_0(\x) p(\x,0) \dd \x,
\end{equation}
where  \(\phi_0\) is the initial data of the dispersion model~\labelcref{eq: dispersion model}, \( \xi \) is computed based on the emission model~\labelcref{eq: emission model}, and the adjoint \( p \) is the weak solution to the following adjoint equation:
\begin{subequations}\label{eq: adjoint}
	\begin{align}
		-\partial_t p(\x,t) - \mu \Delta p(\x,t) - \vv \cdot \nabla p(\x,t) 
		+\kappa \p(\x,t)
		& = \frac{1}{T \vert\Omega\vert} 
		&& \text{for } (\x,t) \in \Omega \times (0,T), \\
		\mu \nabla p(\x,t) \cdot \outernormal(\x) 
		+ \vv \cdot \outernormal(\x) \, p(\x,t) & = 0 
		&& \text{for } (\x,t) \in S^+(\vv), \\
		\nabla p(\x,t) \cdot \outernormal(\x) & = 0 && \text{for } (\x,t) \in S^-(\vv), \\
		p(\x,T) & = 0 && \text{for } \x \in \Omega. 
	\end{align}
\end{subequations}
Again, the existence and uniqueness of a solution to the adjoint equation~\labelcref{eq: adjoint} follows from~\cite[Theorem 2]{Skiba2000}, making the adjoint equation well-posed.

\subsection{Multi-objective optimization involving traffic emission models}\label{sec: multi-objective optimization}
After having introduced functions to measure the economic efficiency and environmental impact of vehicular traffic, we now formulate an optimization problem that optimizes both aspects simultaneously. 
In other words, we try to optimize multiple objectives at the same time which leads us to the concept of multi-objective optimization.

The key components of a multi-objective optimization problem are the same as those of a classic, single-objective one:
\emph{Controls} \( \V \in \R^d \), a \emph{feasible set} \( \S \subseteq \R^d \) encoding the constraints of the optimization if there are any, 
and \( m \ge 2 \) \emph{objectives} \( \J_i \colon \S \rightarrow \R \) to be minimized. 
The main difference to classical optimization is that we have multiple scalar-valued objectives \( \J_i \) instead of a single one and we want to minimize all of them at once.
The corresponding multi-objective optimization problem reads:
\begin{equation}\label{eq: optimization framework}
	\min_{\V \in \S}
	\JJ(\V) \coloneqq \left( \J_1(\V), \ldots, \J_m(\V) \right).
\end{equation}
The objectives \( \J_i\) are typically conflicting which means that there exists no control \( \tilde{\V}\) such that \( \J_i(\tilde{\V}) = \min_{\V \in \S} \J_i(\V)\) for all \( i = 1, \ldots, m\). 
Hence, we characterize solutions to \cref{eq: optimization framework} with the help of Pareto optimality. For a comprehensive overview on multi-objective optimization and the concept of Pareto optimality see \cite{Branke2008,Ehrgott2005} and the references therein.

A control \( \V^{\ast}\) is called \emph{Pareto optimal} if no component of \( \JJ(\V^{\ast})\) can be improved without deteriorating at least one of the others. 
The corresponding objective value \( \JJ(\V^{\ast})\) is called \emph{efficient} and the set of all efficient points is called \emph{Pareto front}.
Thus, solving a multi-objective optimization problem correspond to the following task:
\begin{equation*}%\label{eq: pareto optimal}
	\min_{\V \in \S} \, \JJ(\V) \quad \logeq \quad
	\begin{cases} \text{Find all Pareto optimal controls } \V \in \S, & \\
		\text{and compute the Pareto front of } \JJ. & \end{cases} 
\end{equation*}

\medbreak
Now, we can formulate suitable multi-objective optimization problems that reflect our set goal of optimizing economic and environmental aspects of vehicular traffic at the same time by adjusting the speed limits on the roads. 
Our first ansatz is to measure economic efficiency by the accumulated traffic flow and the environmental impact by the average mass of air pollutant emitted by active and idle traffic.
Altogether, this yields the following multi-objective optimization problem:
\begin{subequations}\label{eq: optimization problem}
	\begin{gather}
	\min_{\V^{\max} \in \S} \, \left( -\J_\flow(\V^{\max}), \, \J_\poll(\V^{\max}; \delta) \right) \\
	\text{where} \quad \S = \Big \{ \V^{\max}  \in (0, \infty)^d \mathrel{\big|}
		 \underline{V_e} \le V_e^{\max} \le \overline{V}_e \Big \}. \label{eq: constraints}
	\end{gather}
\end{subequations}
The scalar \( d\) denotes the number of roads within the road network. We also assume that the speed limit on each road is bounded from below and above by \( \underline{V_e} > 0\) and \( \overline{V}_e < \infty\), respectively to avoid \( V^{\max}_e = 0\) and to mirror that traffic cannot be infinitely fast.
Another ansatz is to minimize the environmental impact of active and idle traffic individually, i.e., splitting the objective \( \J_\poll\) into two, resulting in the multi-objective optimization problem:
\begin{subequations}\label{eq: optimization problem 3D}
	\begin{gather}
		\min_{\V^{\max} \in \S} \, \left( -\J_\flow(\V^{\max}), \, \J_\diff(\V^{\max}), \J_\queue(\V^{\max})\right) \\
		\text{where} \quad \S = \Big \{ \V^{\max}  \in (0, \infty)^d \mathrel{\big|}
		\underline{V_e} \le V_e^{\max} \le \overline{V}_e \Big \},
	\end{gather}
\end{subequations}
with three objectives instead of two compared to the problem in \labelcref{eq: optimization problem}.

The remainder of this paper focuses on solving the proposed multi-objective optimization problems, i.e., the computation of their Pareto fronts. We approach this in a \emph{first-discretize-then-optimize} manner.% and hence, our next concern is to discretize the traffic, the emission and the dispersion model as well as the objectives.
\section{Numerical discretization}\label{sec: numerics}
In this section, we derive discrete counterparts to our optimization problems~\labelcref{eq: optimization problem,eq: optimization problem 3D}. This requires the discretization of the traffic model~\labelcref{eq: lwr}, the emission model~\labelcref{eq: emission model}, and the adjoint equation~\labelcref{eq: adjoint} (as we invoke the expression~\labelcref{eq: alternative expression J_diff} for \( \J_\diff\)).
With the discretization of the above, we can derive the discrete formulations for the objectives \( \J_\flow \), \(\J_\diff\), \(\J_\queue\), and \( \J_\poll\).

Because the emission model depends on the solution to the traffic model and the spatial domain of the adjoint equation, we cover its discretization after having introduced the numerical treatment of the traffic model and the adjoint equation.
Further, we omit the dependence on the speed limit policy \( \V^{\max} \) in the following to improve the readability.

%The evaluation of the objectives requires the explicit solution of the traffic and dispersion model -- or the adjoint equation depending on the approach to compute the function \( \J_\diff \).
%In either case, the involved PDEs lack a closed-form solution in general, and therefore, we opt to solve them numerically.

\subsection{The traffic model}
Before we discretize the traffic model with a finite volume scheme, cf. \cite{LeVeque1992}, we require a discretization of its temporal and spatial domain first.  
Therefore, we approximate the time horizon \( [0,T]\) by an equidistant grid with \( N_t \) grid points and denote them by
\begin{equation}\label{eq: time grid}
	t^k \coloneqq k \Delta t \quad \text{with} \quad \Delta t \coloneqq \frac{T}{N_t},
\end{equation}
where \( \Delta t \) is called \emph{temporal step size}. 
For the discretization of the roads, we assume that all roads have the same length, i.e., the interval corresponding to the road \( e \) is given by \( I_e = I \) where \( I = [0, L^\lwr] \).
Then, we divide \( I \) into \( N_s \) cells \( I_n \) of equal length \( \Delta s \) which means
\begin{equation}\label{eq: grid lwr}
	I_n \coloneqq [s_{n-1}, s_n) \quad \text{where} \quad
		s_n = n \Delta s, \quad \text{and}  \quad
		\Delta s = \frac{L^\lwr}{N_s}.
\end{equation}

Now, we can use a first-order finite volume scheme on each road to approximate the cell averages of the traffic densities on the cells \(I_n\). For road \(e\), we approximate the cell average of cell \( I_n\) at time \( t^k\) by \( \rho_{e,n}^k \), i.e.,
\begin{equation*}
	\rho_{e,n}^k \approx \frac{1}{\Delta s} \int_{I_n} \rho_e(s,t^k) \dd s.
\end{equation*}
Under the additional assumption that the real solution \( \rho_e\) is constant on each cell, we obtain the relation \( \rho_e(s,t^k) \approx \rho_{e,n}^k\) if \( s \in I_n\).
Then, a finite volume scheme applied to the conservation laws of the traffic model~\labelcref{eq: lwr} reads:
\begin{subequations}\label{eq: FV boundary}
	\begin{align}
		\rho_{e,1}^{k+1} & = \rho_{e,1}^k - \frac{\Delta t}{\Delta s} 
		\bigg [ Q_e^\numeric(\rho_{e,1}^k, \rho_{e,2}^k) - Q_e^{\out,k}\bigg ], \\	
		\rho_{e,n}^{k+1} & = \rho_{e,n}^k  - \frac{\Delta t}{\Delta s} 
		\bigg [ Q_e^\numeric(\rho_{e,n}^k,\rho_{e,n+1}^k)
		- Q_e^\numeric(\rho_{e,n-1}^k,\rho_{e,n}^k) \bigg ] 
		\text{ for } 1 < n < N_s, \\	
		\rho_{e,N_s}^{k+1} & = \rho_{e,N_s}^k - \frac{\Delta t}{\Delta s} 
		\bigg [ Q_e^{\inn,k} - Q_e^\numeric(\rho_{e,N_s-1}^k, \rho_{e,N_s}^k) \bigg],
	\end{align}
\end{subequations}
where \( Q_e^\numeric\) is the so-called \emph{numerical flux function} and estimates the flow between two neighboring cells based on their cell averages. We use the Godunov scheme where the corresponding flux function reads:
\begin{equation}\label{eq: flux function}
	Q_e^\numeric(u,v) = \min\{ D_e(u), S_e(v) \}.
\end{equation}
The quantities \( Q_e^{\out, k} \) and \( Q_e^{\inn, k} \)  describe the in- and outflow rates at the boundary of road \( e\) at time \( t^k \), respectively and in the context of the traffic model are either given through external flow rates or internal coupling conditions.
Notice, that we view the in- and outflow from the view of the connecting junctions, i.e., the outflow from the junction to an outgoing road \( e_j\) is equal to the inflow at the beginning of \(e_j\) and the inflow from an incoming road \( e_i\) to the junction is equal to the outflow at the end of \(e_i\).

\begin{remark}
	The temporal step size \( \Delta t\) of the numerical scheme~\labelcref{eq: FV boundary} has to satisfy a CFL condition for every road \( e\). For the numerical flux function~\labelcref{eq: flux function} the CFL condition reads:
	\begin{equation}\label{eq: CFL lwr}
		\frac{\Delta t}{\Delta s} a_e^{\max} \le 1 
		\quad \text{where} \quad
		a_e^{\max} \coloneqq \max \bigg\{ \vert Q_e^{\prime}(\rho) \vert \mathrel{\Big|}
		\min_{s \in I} \rho_e(s,t^k) \le \rho \le \max_{s \in I} \rho_e(s,t^k) \bigg\}.
	\end{equation}
	Thus, we may need to take additional, intermediate time steps to advance from \( t^k\) to \( t^{k+1}\) where each step satisfies the CFL condition~\labelcref{eq: CFL lwr}.
\end{remark}

We provide a time-dependent inflow rate \( q_e^\inn \) for each access road \(e\) and model the corresponding queue by the ODE~\labelcref{eq: queue ODE}. To solve this ODE, we employ the explicit Euler scheme which results in the following equation:
\begin{align}
	\ell_e^{k+1} = \ell_e^k + \Delta t \left[ q_e^\inn(t^k) 
		- \min \bigg \{ q_e^\inn(t^k) 
		+ \frac{\ell_e^k}{\Delta t}, S_e(\rho_{e,1}^k) \bigg \} \right],
\end{align}
with initial condition \( \ell_e^0 = 0 \). 
To discretize the queue's demand \( D_e^\queue\), cf. \cref{eq: queue demand}, we use \( q_e^\queue(t^k) \approx \ell^k_e / \Delta t\).
The computation of the outflow of the queue \(Q_e^{\out, k}\) to an access road then requires the evaluation of \cref{eq: queue inflow} while using the relations \(\ell_e(t^k) \approx \ell_e^k\) and \( \rho_e(0,t^k) \approx \rho_{e,1}^k\).
In case \( e\) is an exit road, we obtain \(Q_e^{\inn,k}\) by recalling that we assume free flow at the end of the exist roads and utilizing the relation \(\rho_e(L^\lwr,t^k) \approx \rho_{e,N_s}^k\).
Similarly, we compute the in- and outflow rates at internal junctions where we need to evaluate the equations~\labelcref{eq: 1to1 junction}, \labelcref{eq: 1to2 junction,eq: 2to1 junction}.

\subsection{The formal discrete adjoint equation}
We require the solution to the adjoint equation~\labelcref{eq: adjoint} to evaluate \( \J_\diff\) according to \cref{eq: alternative expression J_diff}. In order to utilize standard numerical methods for advection-diffusion-type equations, we transform \( t \mapsto T - t\) to obtain the time-reversed adjoint equation which is formulated forward in time and which we solve by using finite differences.

We exemplify the numerical treatment by restricting the spatial domain to squares, i.e., for a fixed length \( L \) we have \( \Omega = [0,L] \times [0,L]\) and by assuming that the diffusion coefficient \( \mu\) and the wind field \( \vv\) are constant over time.
This implies that the in- and outflow boundaries of \( \Omega \) remain unchanged. 
Further, we denote the time-reversed adjoint by \( \tilde{p}\).

\medbreak
First, we discretize the temporal and spatial domain, before delving into the numerical scheme in more detail.
Since the spatial domain \( \Omega \) is two-dimensional, we assume the number of grid points is the same in both directions. Therefore, we fix a number \( N_h \) of grid points and introduce the notation
\begin{equation}\label{eq: grid omega}
	x_i \coloneqq i h, \quad y_j \coloneqq j h \quad \text{with} \quad h \coloneqq \frac{L}{N_h},
\end{equation}
where \( h \) is called \emph{spatial step size}.
This leads to the discrete and equidistant spatial domain \( \Omega_h \) which is defined as follows
\begin{align*}
	\Omega_h \coloneqq  \{ (x_i, y_j) \mid 0 \le i,j \le N_h \}.
\end{align*}
The temporal domain \( [0,T]\) is discretized in the same fashion as for the traffic model, cf. \cref{eq: time grid}.

\subsubsection{Discretization of the time-reversed adjoint equation}
The idea of finite differences is to approximate the solution of an PDE at the grid points.
Hence, we introduce the notation \( \p(x_i,y_j,t^k) \approx \p_{i,j}^k \) where \( \p_{i,j}^k \) is the numerical solution of the time-reversed adjoint equation which approximates the actual solution \( \p \). An approximation of the adjoint state \( p\) is then given through the relation:
\begin{equation*}
	p_{i,j}^k \approx p(x_i,y_j,t^k) = \p(x_i,y_j,T-t^k) \approx \p_{i,j}^{N_t - k}.
\end{equation*}

We discretize the time-reversed adjoint equation componentwise by approximating the temporal derivative with a forward difference:
\begin{equation*}
	\partial_t \p(x_i,y_j,t^k) 
	\approx \frac{1}{\Delta t} \left( \p_{i,j}^{k+1} - \p_{i,j}^k \right),
\end{equation*}
the Laplace operator \( \mu \Delta \p \) with the standard five-point stencil:
\begin{equation*}
	\mu \Delta \p(x_i,y_j,t^k) 
		\approx	\Delta_\diff \p_{i,j}^k  
		\coloneqq \frac{\mu}{h^2} \left[ \p_{i,j+1}^k + \p_{i,j-1}^k
			-4 \p_{i,j}^k + \p_{i+1,j}^k + \p_{i-1,j}^k \right],
\end{equation*}
and the advection term \( -\vv \cdot \nabla \p \), with the upwind scheme for both, the derivative with respect to \( x \) and with respect to \( y \).
We apply the scheme to \( \tilde{\vv}  \cdot \p\), where \( \tilde{\vv} = - \vv\) to depict the correct direction of advection in the discrete equation, i.e.,
%cf. \cite{Numabibel}: 
\begin{equation*}
	\tilde{\vv} \cdot \nabla \p(x_i,y_j, t^k)
	\approx \Delta_\adv \p_{i,j}^k 
	\coloneqq \frac{1}{h} \left[ \tilde{v}_y^- \p_{i,j+1}^k - \tilde{v}_y^+ \p_{i,j-1}^k 
	+ \Vert \tilde{\vv} \Vert_1 \p_{i,j}^k 
	+ \tilde{v}_x^- \p_{i+1,j}^k - \tilde{v}_x^+ \p_{i-1,j}^k \right],
\end{equation*}
where for a scalar \( a\), the identities \( a^+ \coloneqq \max \{ a, 0 \} \) and \( a^- \coloneqq \min \{a, 0 \} \) hold.
Altogether, by combining the discretization of these different components, we obtain an equation to compute the numerical solution for the next time step \( t^{k+1} \) provided the numerical solution at the current time step \( t^k \) is known:
\begin{equation}\label{eq: adv diff numerical scheme}
	\frac{1}{\Delta t} \left(\p_{i,j}^{k+1} - \p_{i,j}^k \right) 
	- \Delta_\diff \p_{i,j}^{k} + \Delta_\adv \p_{i,j}^k + \kappa \p_{i,j}^k  = \frac{1}{T \vert \Omega \vert}.
\end{equation}
%Note that the numerical solution \( \p_{i,j}^0\) is known as we impose initial data.

\begin{remark}\label{remark cfl}
	Under the assumption that \( \kappa = 0\), and that the wind field \(\vv = (v_x, v_y)^\top\) and the diffusion coefficient \( \mu \) in~\labelcref{eq: dispersion model} are constant over time, the CFL condition corresponding to the discretization~\labelcref{eq: adv diff numerical scheme} follows directly from \cite[Theorem 3.1 and Section 4]{Wesseling1996} and reads:
	\begin{equation*}
		\Delta t \le \frac{h^2}{4 \mu + \Vert \vv \Vert_1 h} \quad \text{and} \quad
		\Delta t \left[ \frac{v_x^2}{2\mu + \vert v_x \vert h} 
		+ \frac{v_y^2}{2 \mu + \vert v_y \vert h}\right] \le 1.
	\end{equation*} 
	However, this CFL condition is invalid for \( \kappa \neq 0\). Hence, we tighten the CFL condition to:
	\begin{equation}\label{eq: cfl adjoint}
		\Delta t \le \frac{1}{3} \frac{h^2}{4 \mu + \Vert \vv \Vert_1 h} \quad \text{and} \quad
		\Delta t \left[ \frac{v_x^2}{2\mu + \vert v_x \vert h} 
		+ \frac{v_y^2}{2 \mu + \vert v_y \vert h}\right] \le \frac{1}{3}.
	\end{equation}
	We found that for our numerical experiments, the discretization~\labelcref{eq: adv diff numerical scheme} is stable given the tightened CFL condition~\labelcref{eq: cfl adjoint}.
	As this condition is independent of \( \kappa \), one should not expect that it holds for all \( \kappa\).
\end{remark}

\subsubsection{Discretization of the boundary conditions}
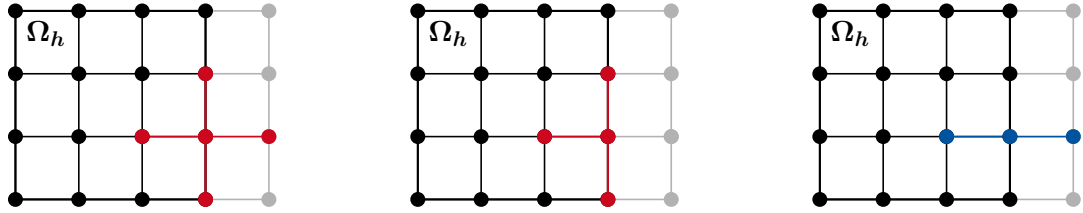
\begin{figure}[ht]
	\centering
	\begin{subfigure}[t]{.32\textwidth}
		\centering
		\scalebox{1}{
	\begin{tikzpicture}[
		semithick,
		]
		
		\pgfmathsetmacro{\N}{3}
		\pgfmathsetmacro{\L}{2.5}
		\pgfmathsetmacro{\h}{\L / \N}

		% spacial grid
		\draw[black] (0,0) grid[step=\h] (\L, \L);
		\draw[black, thick] (0,0) rectangle (\L, \L);

		% road network
		\tikzstyle{state}=[
		draw = black,
		fill=black, 
		circle,
		inner sep=0pt,
		outer sep=0pt,
		minimum size=5pt
		]
		
		\tikzstyle{ghost point}=[
		draw = black!30!white,
		fill = black!30!white, 
		circle,
		inner sep=0pt,
		outer sep=0pt,
		minimum size=5pt
		]
		
		\tikzstyle{stencil}=[
		draw = rot,
		fill = rot, 
		circle,
		inner sep=0pt,
		outer sep=0pt,
		minimum size=5pt
		]
		
			% ghost points in grey
		\foreach \i in {0,...,\N}{
			% east
			\draw[black!30!white] ( \L + \h, \i * \h) -- (\L, \i * \h);
			\node[ghost point] at ( \L + \h, \i * \h) {};
		}
		
		\path[-, black!30!white] (\L + \h, 0) edge node[midway, right] {} (\L + \h, \L);
		
		% colored line for the stencil
		\foreach \i in {-1,0}{
			\draw[rot] (\L + \i * \h, \L - 2 * \h) -- ({\L + (\i + 1) * \h}, \L - 2 * \h);
		}
		\foreach \j in {-1,-2}{
			\draw[rot, thick] (\L , \L + \j * \h) -- (\L, {\L + (\j - 1) * \h});
		}

		% dots at inner grid
		\foreach \i in {0,...,\N}{
			\foreach \j in {0,...,\N}{
				\node[state] at (\i * \h, \j * \h) {};
			}
		}

		% stencil at the left boundary to show inadmissablility
		\foreach \i in {-1,0,1}{
			\node[stencil] at (\L + \i * \h, \L - 2 * \h) {};
		}
		
		\foreach \j in {-1,-2,-3}{
			\node[stencil] at (\L , \L + \j * \h) {};
		}
		
		 \node[below right] at (0, \L) {$\bm{\Omega_h}$};
		
	\end{tikzpicture}
}
		\subcaption{A general five-point stencil.}
		\label{fig: grid ghost cells a}
	\end{subfigure} \hfill
	\begin{subfigure}[t]{.32\textwidth}
		\centering
		\scalebox{1}{
	\begin{tikzpicture}[
		semithick,
		]
		
		\pgfmathsetmacro{\N}{3}
		\pgfmathsetmacro{\L}{2.5}
		\pgfmathsetmacro{\h}{\L / \N}

		% spacial grid
		\draw[black] (0,0) grid[step=\h] (\L, \L);
		\draw[black, thick] (0,0) rectangle (\L, \L);
		
		% road network
		\tikzstyle{state}=[
		draw = black,
		fill=black, 
		circle,
		inner sep=0pt,
		outer sep=0pt,
		minimum size=5pt
		]
		
		\tikzstyle{ghost point}=[
		draw = black!30!white,
		fill = black!30!white, 
		circle,
		inner sep=0pt,
		outer sep=0pt,
		minimum size=5pt
		]
		
		\tikzstyle{stencil}=[
		draw = rot,
		fill = rot, 
		circle,
		inner sep=0pt,
		outer sep=0pt,
		minimum size=5pt
		]
		
		\path[-, black!30!white] (\L + \h, 0) edge node[midway, right] {} (\L + \h, \L);
				
		% ghost points in grey
		\foreach \i in {0,...,\N}{
			% east
			\draw[black!30!white] ( \L + \h, \i * \h) -- (\L, \i * \h);
			\node[ghost point] at ( \L + \h, \i * \h) {};
		}
		
		% coloered line for the stencil
		\foreach \i in {-1}{
			\draw[rot] (\L + \i * \h, \L - 2 * \h) -- ({\L + (\i + 1) * \h}, \L - 2 * \h);
		}
		\foreach \j in {-1,-2}{
			\draw[rot, thick] (\L , \L + \j * \h) -- (\L, {\L + (\j - 1) * \h});
		}
		
		% ghost dot
		\node[ghost point] at ( \L + \h, \L - 2 * \h) {};
		
		% dots at inner grid
		\foreach \i in {0,...,\N}{
			\foreach \j in {0,...,\N}{
				\node[state] at (\i * \h, \j * \h) {};
			}
		} 
		
		% stencil at the left boundary to show inadmissablility
		\foreach \i in {-1,0}{
			\node[stencil] at (\L + \i * \h, \L - 2 * \h) {};
		}
		
		\foreach \j in {-1,-2,-3}{
			\node[stencil] at (\L , \L + \j * \h) {};
		}
		
		 \node[below right] at (0, \L) {$\bm{\Omega_h}$};

	\end{tikzpicture}
}
		\subcaption{A reduced five-point stencil.}
		\label{fig: grid ghost cells c}
	\end{subfigure} \hfill
	\begin{subfigure}[t]{.32\textwidth}
		\centering
		\scalebox{1}{
	\begin{tikzpicture}[
		semithick,
		]
		
		\pgfmathsetmacro{\N}{3}
		\pgfmathsetmacro{\L}{2.5}
		\pgfmathsetmacro{\h}{\L / \N}

		% spacial grid
		\draw[black] (0,0) grid[step=\h] (\L, \L);
		\draw[black, thick] (0,0) rectangle (\L, \L);
		
		% road network
		\tikzstyle{state}=[
		draw = black,
		fill=black, 
		circle,
		inner sep=0pt,
		outer sep=0pt,
		minimum size=5pt
		]
		
		\tikzstyle{ghost point}=[
		draw = black!30!white,
		fill = black!30!white, 
		circle,
		inner sep=0pt,
		outer sep=0pt,
		minimum size=5pt
		]
		
		\tikzstyle{stencil}=[
		draw = rwth,
		fill = rwth, 
		circle,
		inner sep=0pt,
		outer sep=0pt,
		minimum size=5pt
		]
		
		\path[-, black!30!white] (\L + \h, 0) edge node[midway, right] {} (\L + \h, \L);
		
		% ghost points in grey
		\foreach \i in {0,...,\N}{
			% east
			\draw[black!30!white] ( \L + \h, \i * \h) -- (\L, \i * \h);
			\node[ghost point] at ( \L + \h, \i * \h) {};
		}
		
		% coloered line for the stencil
		\foreach \i in {-1,0}{
			\draw[rwth] (\L + \i * \h, \L - 2 * \h) -- ({\L + (\i + 1) * \h}, \L - 2 * \h);
		}
		
		% ghost dot
		\node[ghost point] at ( \L + \h, \L - 2 * \h) {};
		
		% dots at inner grid
		\foreach \i in {0,...,\N}{
			\foreach \j in {0,...,\N}{
				\node[state] at (\i * \h, \j * \h) {};
			}
		} 
		
		% stencil at the left boundary to show inadmissablility
		\foreach \i in {-1,0,1}{
			\node[stencil] at (\L + \i * \h, \L - 2 * \h) {};
		}	
		
		  \node[below right] at (0, \L) {$\bm{\Omega_h}$};
	\end{tikzpicture}
}
		\subcaption{The central-difference stencil.}
		\label{fig: grid ghost cells b}
	\end{subfigure}
	\caption{%Visualization of the procedure to handle possible ghost points when applying a five-point stencil to the right boundary of a square domain.
	The black dots represent the grid, the red and blue dots a (reduced) five-point stencil and the stencil of a central difference respectively. The grey dots indicate the ghost points.}
	\label{fig: grid ghost cells}
\end{figure}

The numerical scheme~\labelcref{eq: adv diff numerical scheme} suffers from a shortcoming: the stencils \( \Delta_\diff \) and \( \Delta_\adv \) applied to boundary values demand values of the numerical solution at inadmissible points, so-called \emph{ghost points}, outside the domain \( \Omega \) as illustrated in \cref{fig: grid ghost cells a}. 
%These points are so-called \emph{ghost points} and the values associated with the numerical solution at those points are called \emph{ghost values}.
In the following, we provide a brief description of how to resolve this problem.

The (time-reversed) adjoint equation involves two types of boundary conditions: Neumann boundary conditions and Robin boundary conditions.
Both determine the value of the solution implicitly through a differential equation.

The key idea to resolve the problem, is to utilize the ghost points when discretizing the Neumann and Robin boundary condition. We illustrate the method for the case that the right boundary of \( \Omega\) is a Robin boundary and discretize the corresponding Robin boundary condition as follows:
\begin{equation*}
	\mu \underbrace{\frac{\p_{N_h+1,j}^k - \p_{N_h-1,j}^k}{2 h}}_{\approx \partial_x \p(x_{N_h}, y_j,t^k)} 
	+ v_x \underbrace{\frac{\p_{N_h+1,j}^k + \p_{N_h-1,j}^k}{2}}_{\approx \p(x_{N_h},y_j,t^k)}
	= 0.
\end{equation*}
The corresponding stencil is visualized in \cref{fig: grid ghost cells b}. Now, we can express the value \( \p_{N_h+1,j}^k\) at a ghost point by an admissible value of the numerical solution:
\begin{equation}\label{eq: ghost value}
	\p_{N_h+1,j}^k 
		= \frac{\mu - v_x h}{\mu + v_x h} \p_{N_h-1,j}^k.
\end{equation}
Finally, we can use the relation~\labelcref{eq: ghost value} to eliminate the inadmissable values from the stencils \( \Delta_\diff \p_{N_h,j}^k\) and \( \Delta_\adv \p_{N_h,j}^k\) which results in a reduced five-point stencil as illustrated in \cref{fig: grid ghost cells c}.
The approach in case the right boundary is a Neumann boundary is analogous and yields:
\begin{equation*}%\label{eq: ghost value}
	\p_{N_h+1,j}^k 	= \p_{N_h-1,j}^k.
\end{equation*}
 %Then, we can express the ghost value by admissible values of the numerical solution and can use the previously established relation to eliminate the ghost value from the stencil. 
%\cref{fig: grid ghost cells c} shows the reduced stencil one obtains from this ansatz.

%
\subsection{The emission model}\label{sec: numerics emission}
After having established the numerical treatment of the traffic model~\labelcref{eq: lwr} and the adjoint equation~\labelcref{eq: adjoint}, we proceed with approximating the emission rate \( \xi\) on \( \Omega_h\) at given times \( t^k\).
Under the additional assumption that the traffic densities \( \rho_e \) are piecewise constant on the cells \( I_n \), we approximate the emission rate of road \( e \)  at time \( t^k \) by
\begin{align*}
	\xi_e(s,t^k) \approx \tilde{\xi}_{e,n}^k  \coloneqq Q_e(\rho_{e,n}^k) + \theta \rho_{e,n}^k
	\quad \text{if} \quad s \in I_n.
\end{align*}
Given these approximations, we invoke \cref{eq: emission model}, which directly yields an approximation of the emission rate \( \xi \) at a point \( \x \in \Omega \) where \( \E(\x) \neq \emptyset\) holds:
\begin{equation*}
	\xi(\x,t^k) \approx \xi^h(\x,t^k) \coloneqq \frac{1}{\vert \E(\x) \vert} 
	\sum_{e \in \E(\x)} \frac{1}{w_e} \tilde{\xi}_{e,n}^k
	\quad \text{if } \x \in R_e(s) \text{ and } s \in I_n.
\end{equation*}
Then, at a grid point \( (x_i,y_j)\) of the discrete control area, the approximation of the emission rate reads:
\begin{equation}\label{eq: discrete emission model}
	\xi(x_i,y_j,t^k) \approx \xi_{i,j}^k \coloneqq 
	\begin{cases}	
		\xi^h(x_i,y_j,t^k), & \text{if } \E(x_i,y_j) \neq \emptyset,  \\
		0, & \text{otherwise.}
	\end{cases}
\end{equation}

\begin{figure}[ht]
	\centering
	\begin{subfigure}[t]{.475\textwidth}
		\centering
		\scalebox{1.0}{
\begin{tikzpicture}[
	> = {Latex[scale=1.0]},
    thick,
]
        
\pgfmathsetmacro{\Nlwr}{5}
\pgfmathsetmacro{\Llwr}{4}
\pgfmathsetmacro{\hlwr}{\Llwr / \Nlwr}
\pgfmathsetmacro{\L}{\Llwr}
\pgfmathsetmacro{\N}{3}
\pgfmathsetmacro{\h}{\L / \N}
\pgfmathsetmacro{\xx}{\Llwr + 6}
\pgfmathsetmacro{\w}{\L/10}
\pgfmathsetmacro{\hh}{1.5}

% road network       
\tikzstyle{state}=[
   draw= black,
   fill = black, 
   circle,
   inner sep=0pt,
   outer sep=0pt,
   minimum size=5pt
 ]
 
\tikzstyle{grey state}=[
	draw = black!50!white,
    fill = black!50!white, 
    circle,
    inner sep=0pt,
    outer sep=0pt,
    minimum size=5pt
]

	\draw[fill=black!10!white, draw=black!10!white, thick] (\L/2-\w,0) rectangle (\L/2 + \w,\L);
	
	\draw[black!100!white] (0,0) grid[step=\h] (\L, \L);
	\draw[black!100!white, thick] (0,0) rectangle (\L, \L);
		   
	\path[->, black!50!white, thick, shorten >= 3pt] (0.5*\L, 0) edge node[midway, right] {} (0.5*\L, \L);  
	\foreach \i in {0,...,\N}{
		\foreach \j in {0,...,\N}{
			\node[state] at (\i * \h, \j * \h) {};
		}
	}   
	
	% code color the grid ponits but in 2d
	\foreach \i in {0,...,\Nlwr}{
		\node[grey state] at (\Llwr / 2, \i * \hlwr) {};
	} 
	\draw[decoration={brace, mirror, raise=5pt, amplitude=5pt}, decorate]%,line width=1pt]	
  		(\L / 2, 0) -- (\L / 2, \hlwr) node[midway, right=10pt]{$\Delta s$};
	
	\draw[decoration={brace, mirror, raise=5pt, amplitude=5pt}, decorate]	
  		(\L,0) -- (\L,\h) node[midway, right=10pt]{$h$};
  		
  	\node[below right] at (0, \L) {$\bm{\Omega_h}$};
  			
\end{tikzpicture}
}
		\subcaption{Unsuitable spatial step sizes.} 
	\end{subfigure} \hfill
	\begin{subfigure}[t]{.475\textwidth}
		\centering
		\scalebox{1.0}{
\begin{tikzpicture}[
	> = {Latex[scale=1.1]},
    thick,
    %shorten >= 3pt
]
        
\pgfmathsetmacro{\Nlwr}{5}
\pgfmathsetmacro{\Llwr}{4}
\pgfmathsetmacro{\hlwr}{\Llwr / \Nlwr}
\pgfmathsetmacro{\L}{\Llwr}
\pgfmathsetmacro{\N}{5}
\pgfmathsetmacro{\h}{\L / \N}
\pgfmathsetmacro{\xx}{\Llwr + 6}
\pgfmathsetmacro{\w}{\L/10}
\pgfmathsetmacro{\hh}{1.5}

% road network       
\tikzstyle{state}=[
   draw= black,
   fill = black, 
   circle,
   inner sep=0pt,
   outer sep=0pt,
   minimum size=5pt
 ]
 
\tikzstyle{grey state}=[
	draw = black!50!white,
    fill = black!50!white, 
    circle,
    inner sep=0pt,
    outer sep=0pt,
    minimum size=5pt
]
		
	\draw[fill=black!10!white, draw=black!10!white, thick] (\L/2-\w,0) rectangle (\L/2 + \w,\L);
	
	\draw[black!100!white] (0,0) grid[step=\h] (\L, \L);
	\draw[black!100!white, thick] (0,0) rectangle (\L, \L);
		   
	\path[->, black!50!white, thick,shorten >= 3pt] (0.5*\L, 0) edge node[midway, right] {} (0.5*\L, \L);  
	\foreach \i in {0,...,\N}{
		\foreach \j in {0,...,\N}{
			\node[state] at (\i * \h, \j * \h) {};
		}
	}   
	
	% code color the grid ponits but in 2d
	\foreach \i in {0,...,\Nlwr}{
		\node[grey state] at (\Llwr / 2, \i * \hlwr) {};
	} 
	\draw[decoration={brace, mirror, raise=5pt, amplitude=5pt}, decorate]%,line width=1pt]	
  		(\L / 2, 0) -- (\L / 2, \hlwr) node[midway, right=10pt]{$\Delta s$};
	
	\draw[decoration={brace, mirror, raise=5pt, amplitude=5pt}, decorate]%,line width=1pt]	
  		(\L, 0) -- (\L, \h) node[midway, right=10pt]{$h$};
  		
  	  	\node[below right] at (0, \L) {$\bm{\Omega_h}$};

\end{tikzpicture}
}
		\subcaption{Suitable spatial step sizes.}
		%\label{fig: suitable step size}
	\end{subfigure}
	\caption{Illustration of suitable and unsuitable discretizations of \(\Omega\) for a single-roaded network.}
	\label{fig: unsuitable step size}
\end{figure}
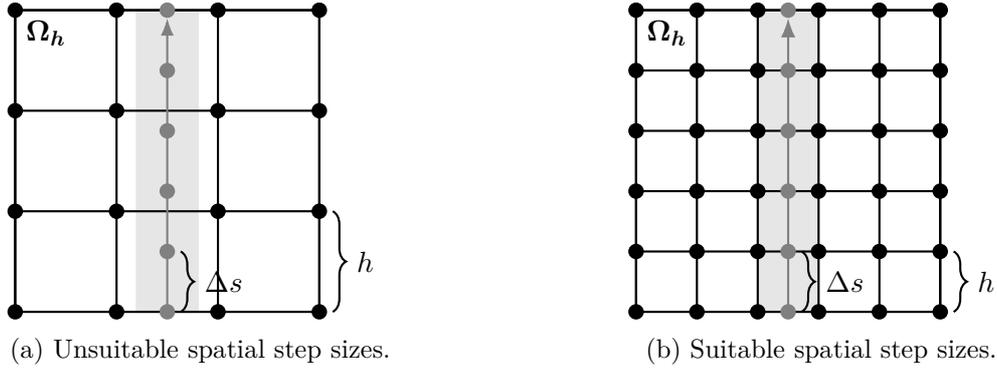

\begin{remark}
	The spatial step size of the advection-diffusion-type equation cannot be set arbitrarily. If, for example, the step size is too large, the road network is invisible to the discrete control area. 
	We visualize the issue for a road network consisting of a single road in \cref{fig: unsuitable step size} including a discretization that is consistent with the length and width of the road as well as its positioning. 
\end{remark}

%To close the discussion on the discretization of the traffic emission model and the adjoint equation, we summarize the order of computations to obtain numerical solutions to the model variables and the adjoint state:
%\begin{enumerate}[label=(\roman*)]
%	\item Set the step sizes \( \Delta s\), \( h \), and \( \Delta \tau \). 
%	\item Solve the traffic model with \cref{alg: LWR} and obtain the traffic densities 
%	\( \rho_{e,n}^k \) and queue length' \( \ell_e^k \). \;
%	\item Compute the emission rate with \cref{eq: discrete emission model}.
%	\item Solve the advection-diffusion-type equation with \cref{alg: advection diffusion} to either obtain
%	the concentration level \( \phi_{i,j}^k \), or adjoint state \( p_{i,j}^k \). 
%\end{enumerate}
%Since the adjoint equation is independent of the speed limits and the model variables, its numerical solution can also be computed before solving the traffic model.
%We utilize these numerical solutions to discretize the objectives with a quadrature rule in the next section.

\subsection{The objectives and the discrete optimization problem}\label{sec: numerics objectives}
Now given  numerical solutions of the traffic model~\labelcref{eq: lwr}, the emission model~\labelcref{eq: emission model}, and the adjoint equation~\labelcref{eq: adjoint} with the schemes introduced above, we derive discrete versions of the objectives \( \J_\flow\), \(\J_\diff\), \( \J_\queue\), and \( \J_\poll\) based on the extended right rectangular rule for integrals.

Recall that the continuous variables depend on the speed limits, except for the adjoint state. Thus, the discrete variables inherit this dependence, which means we can write \( \rho_{e,n}^k = \rho_{e,n}^k(\V^{\max}) \), \( \ell_e^k = \ell_e^k(\V^{\max}) \), and \( \xi_{i,j}^k = \xi_{i,j}^k(\V^{\max}) \).
This leads to the following discrete functions indicated by the superscript \( h\):
%
%\medbreak
%First, we consider the economic objective that measures the accumulated traffic flow. Since we assume  that the traffic densities are piecewise constant on the cells \( I_n\), we obtain:
\begin{subequations}
\begin{align}
	\J_\flow^h(\V^{\max})& =
		\Delta t h \sum_{e \in E} \sum_{k=1}^{N_t} \sum_{n=1}^{N_s} Q_e(\rho_{e,n}^k(\V^{\max}), \V^{\max}), \\
	\J_\diff^h(\V^{\max}) & = \Delta t h^2 	\sum_{k=1}^{N_t} \sum_{i,j=1}^{N_h} \xi_{i,j}^k(\V^{\max}) p_{i,j}^k
	+ h^2 \sum_{i,j=1}^{N_h} \phi_0(x_i,y_j) p_{i,j}^0 \label{eq: J_diff discrete} \\
	\J_\queue^h(\V^{\max}) & = \frac{\Delta t}{T} \sum_{e \in E_\inn} \sum_{k=1}^{N_t}  \ell_e^k(\V^{\max}) \\
	\J_\poll^h(\V^{\max}; \delta) & = \J_\diff^h(\V^{\max}) + \delta \J_\queue^h(\V^{\max}). \vphantom{\sum_{k=1}^N} \label{eq: J_poll discrete}
\end{align}
\end{subequations}
%Next, we derive the discrete ecological objective. To do so, we derive discrete expressions of its two components -- the functions measuring the contribution of active and passive traffic to air pollution.
%
%To discretize the function describing the contribution of passive traffic, we get:
%\begin{align*}
%	\J_\queue(\V^{\max};\kappa) & = \frac{\theta}{T } \sum_{e \in E_\inn} \int_0^T \ell_e(t,\V^{\max}) \dd t  \\
%	& \approx \frac{\kappa \Delta t}{T} \sum_{e \in E_\inn} \sum_{k=1}^{N_t}  \ell_e^k(\V^{\max})
%	\eqqcolon \J_\queue^h(\V^{\max}; \kappa),
%\end{align*}
%and for the discrete version describing the contribution of active traffic, the discretization reads:
%\begin{align*}
%	\J_\diff(\V^{\max}) & = \int_0^T \int_{\Omega} \xi(\x,t,\V^{\max}) p(\x,t) \dd \x \dd t
%	+ \int_{\Omega} \phi_0(\x) p(\x,0) \dd \x \\ 
%	& \approx \Delta t h^2 	\sum_{k=1}^{N_t} \sum_{i,j=1}^{N_h} \xi_{i,j}^k(\V^{\max}) p_{i,j}^k
%	+ h^2 \sum_{i,j=1}^{N_h} \phi_0(x_i,y_j) p_{i,j}^0
%	\eqqcolon \J_\diff^h(\V^{\max}).
%\end{align*}
%Combing both discretizations leads to the discrete counterpart of the ecological objective:
%\begin{equation*}
%	\J_\poll^h(\V^{\max};\kappa) = \J_\diff^h(\V^{\max}) + \J_\queue^h(\V^{\max};\kappa).
%\end{equation*}

%\medbreak
Given the discrete counterparts of the objectives, we can finally derive the discrete optimization problems by replacing the continuous objectives with the discrete ones.
The optimization problem~\labelcref{eq: optimization problem}, where the goal is to maximize the traffic flow while simultaneously minimizing the environmental impact, then becomes:
\begin{align}\label{eq: optimization problem discrete}
	\min_{\V^{\max} \in \S} \left( -\J_\flow^h(\V^{\max}), \J_\poll^h(\V^{\max}; \delta)\right),
\end{align}
while the optimization problem~\labelcref{eq: optimization problem 3D} which involves three objectives, becomes:
\begin{align}\label{eq: optimization problem 3D discrete}
	\min_{\V^{\max} \in \S} \left( -\J_\flow^h(\V^{\max}), \J_\diff^h(\V^{\max}), \J_\queue^h(\V^{\max})\right).
\end{align}

The discrete optimization problems~\labelcref{eq: optimization problem discrete,eq: optimization problem 3D discrete} can be solved by standard methods to solve multi-objective optimization problems as the discrete objectives allow a direct evaluation by invoking the numerical schemes introduced previously.

\section{Numerical experiments}\label{sec: experiments}
In this section, we aim to gain an insight into the trade-off between maximizing economic efficiency and minimizing the environmental impact. 
To achieve this, we solve the discretized multi-objective optimization problems~\labelcref{eq: optimization problem discrete,eq: optimization problem 3D discrete} with fixed coefficients and step sizes, providing a proof-of-concept example.
Furthermore, we study the impact of the presence of idle traffic on the solution of the optimization problem~\labelcref{eq: optimization problem}.
We distinguish between two cases: when the function \(\J_\poll \), measuring air pollution, only considers the active traffic in the network (\(\delta = 0\)), and when it also takes the idle traffic in the queues into account (\(\delta > 0\)), cf. \cref{eq: J_poll discrete}.
Additionally, we compare the results to the solution of the optimization problem~\labelcref{eq: optimization problem 3D} which minimizes the contribution of active and idle traffic to air pollution individually.

We solve the resulting multi-objective optimization problems using the \texttt{paretosearch} function from Matlab's \emph{Optimization Toolbox}~\cite{Matlab}.
This function is based on pattern search on a set of points. It iteratively searches for efficient points that are close to the true Pareto front.

\subsection{Definition of the proof-of-concept example}\label{sec: settings}
\begin{figure}[ht]
	\centering
	\begin{subfigure}[t]{.32\textwidth}
		\centering
		\scalebox{0.75}{
	\begin{tikzpicture}[
		> = {Latex[scale=1.1]}, 
		line width=3pt 
		]
		
		\pgfmathsetmacro{\L}{2}
		\pgfmathsetmacro{\h}{0.3}
		
		\tikzstyle{state}=[
		fill=black, 
		circle,
		inner sep=0pt,
		outer sep=0pt,
		minimum size=8pt,
		thick
		]
		
		\tikzstyle{outer state}=[
		draw=black,
		fill=white, 
		circle,
		inner sep=0pt,
		outer sep=0pt,
		minimum size=8pt,
		thick
		]
		
		\draw[black!30!white, fill=black!4!white, line width=1.25pt] (0, -\L) rectangle (3*\L, 2*\L); 
		
		\node[outer state] at (0,0) (1) {};        
		\node[state] at (\L, 0) (2) {}; %{1};
		\node[state] at (2*\L, 0) (3) {}; %{2};
		\node[state] at (2*\L, \L) (4) {}; %{3};
		\node[state] at (\L, \L) (5) {}; %{4};
		\node[outer state] at (3*\L, \L) (6) {};

		\path[->, thick] (1) edge node[midway, below] {$I_1$} (2); %{1} (2);
		\path[->, thick] (2) edge node[midway, below, align=center] 
		{$I_3$} (3); %{3} (3);
		\path[->,  thick] (2) edge node[midway, left, align=center] 
		{$I_2$} (5); %{2} (5);
		\path[->,  thick] (3) edge node[midway, right, align=center]
		{$I_4$} (4); %{4} (4);
		\path[->, thick] (5) edge node[midway, above, align=center] 
		{$I_5$} (4); %{5} (4);
		\path[->,  thick] (4) edge node[midway, above] {$I_6$} (6); %{6} (6);

		\node[outer state] at (0,0) (1) {};        
		\node[state] at (\L, 0) (2) {}; %{1};
		\node[state] at (2*\L, 0) (3) {}; %{2};
		\node[state] at (2*\L, \L) (4) {}; %{3};
		\node[state] at (\L, \L) (5) {}; %{4};
		\node[outer state] at (3*\L, \L) (6) {};
				
		\node[below right] at (0, 2*\L) {$\textcolor{black!50!white}{\bm{\Omega}}$};
	\end{tikzpicture}
}
		\subcaption{Road network in \( \Omega\).}
		\label{fig: road network}
	\end{subfigure} \hfill
	\begin{subfigure}[t]{.32\textwidth}
		\centering
		\scalebox{0.75}{
\begin{tikzpicture}[
	> = {Latex[scale=1.1]}, % arrow head style
    %shorten > = 1pt, % don't touch arrow head to node
    % auto,
    line width=3pt, %ultra thick, % line style
]

\pgfmathsetmacro{\L}{2}
\pgfmathsetmacro{\h}{0.3}

	\tikzstyle{state}=[
    	fill=black, 
        circle,
        inner sep=0pt,
        outer sep=0pt,
        minimum size=8pt,
        thick
   	]
        
     \tikzstyle{outer state}=[
     	draw=black,
     	fill=white, 
        circle,
        inner sep=0pt,
        outer sep=0pt,
        minimum size=8pt,
       	thick
     ]
       
     \draw[black!30!white, fill=black!4!white, line width=1.25pt] (0, -\L) rectangle (3*\L, 2*\L); 

     \node[outer state] at (0,0) (1) {};        
     \node[state] at (\L, 0) (2) {}; %{1};
     \node[state] at (2*\L, 0) (3) {}; %{2};
     \node[state] at (2*\L, \L) (4) {}; %{3};
     \node[state] at (\L, \L) (5) {}; %{4};
     \node[outer state] at (3*\L, \L) (6) {};

     \path[-, draw=black!60!white] (1) edge node[midway, below] {} (2); %{1} (2);
     \path[-, draw=black!90!white] (2) edge node[midway, above, align=center] 
     	{} (3); %{3} (3);
     \path[-, draw=black!40!white] (2) edge node[midway, left, align=center] 
     	{} (5); %{2} (5);
     \path[-, draw=black!50!white] (3) edge node[midway, right, align=center]
     	{} (4); %{4} (4);
     \path[-, draw=black!80!white] (5) edge node[midway, above, align=center] 
     	{} (4); %{5} (4);
     \path[-, draw=black!30!white] (4) edge node[midway, above] {} (6); %{6} (6);
     
  	\path[->, thick] (\h,-\h) edge node[midway, below] {$\rho_1^0 = 0.6$} (\L -\h,-\h);
  	\path[->, thick] (\L - \h,\h) edge node[midway, left] {$\rho_2^0 = 0.4$} (\L - \h, \L - \h);
  	\path[->, thick] (\L + \h,-\h) edge node[midway, below] {$\rho_3^0=0.9$} (2*\L - \h, - \h);
  	\path[->, thick] (2*\L + \h,\h) edge node[midway, right] {$\rho_4^0=0.5$} (2*\L + \h, \L - \h);
  	\path[->, thick] (\L + \h, \L + \h) edge node[midway, above] {$\rho_5^0=0.8$} (2*\L - \h, \L + \h);
  	\path[->, thick] (2*\L + \h,\L+\h) edge node[midway, above] {$\rho_6^0=0.3$} (3*\L - \h, \L + \h);
  	
  	 \node[outer state] at (0,0) (1) {};        
     \node[state] at (\L, 0) (2) {}; %{1};
     \node[state] at (2*\L, 0) (3) {}; %{2};
     \node[state] at (2*\L, \L) (4) {}; %{3};
     \node[state] at (\L, \L) (5) {}; %{4};
     \node[outer state] at (3*\L, \L) (6) {};
     
     \node[below right] at (0, 2*\L) {$\textcolor{black!50!white}{\bm{\Omega}}$};
\end{tikzpicture}
}

%\begin{tikzpicture}[
%	> = {Latex[scale=1.1]}, % arrow head style
 %   %shorten > = 1pt, % don't touch arrow head to node
  %  % auto,
   % node distance = 2.5cm, % distance between nodes
    %semithick, % line style
    %font=\small
%]

%	\tikzstyle{state}=[
 %   	fill=black, 
  %      circle,
   %     inner sep=0pt,
    %    outer sep=0pt,
     %   minimum size=8pt
   	%]
        
     %\tikzstyle{outer state}=[
     %	draw=black, 
      %  circle,
       % inner sep=0pt,
        %outer sep=0pt,
     %   minimum size=8pt
     %]

     %\node[outer state] (1) {};        
     %\node[state] (2) [right of=1] {}; %{1};
     %\node[state] (3) [right of=2] {}; %{2};
     %\node[state] (4) [above of=3] {}; %{3};
     %\node[state] (5) [above of=2] {}; %{4};
     %\node[outer state] (6) [right of=4] {};

     %\path[->] (1) edge node[midway, below] {$\rho_1^0 = 0.6$} (2); %{1} (2);
     %\path[->] (2) edge node[midway, below, align=center] {$\rho_3^0 = 0.9$ \\ $1-\alpha=0.4$} (3); %{3} (3);
    % \path[->] (2) edge node[midway, left, align=center] {$\rho_2^0 = 0.4$ \\ $\alpha = 0.6$} (5); %{2} (5);
   %  \path[->] (3) edge node[midway, right, align=center] {$\rho_4^0 = 0.5$ \\ $1-\beta = 0.4$} (4); %{4} (4);
  %   \path[->] (5) edge node[midway, above, align=center] {$\rho_5^0 = 0.8$ \\$\beta = 0.6$} (4); %{5} (4);
 %    \path[->] (4) edge node[midway, above] {$\rho_6^0 = 0.3$} (6); %{6} (6);
%\end{tikzpicture}
%}
		\subcaption{Initial data.}
		\label{fig: initial data lwr}
	\end{subfigure} \hfill
	\begin{subfigure}[t]{.32\textwidth}
		\centering
		\scalebox{0.75}{
\begin{tikzpicture}[
    thick % line style
]

\pgfmathsetmacro{\L}{2}
\pgfmathsetmacro{\w}{1}
\pgfmathsetmacro{\h}{0.5}
\pgfmathsetmacro{\N}{3*\L / \h}

	\draw[draw=black!30!white, fill=black!4!white] (0,-\L) rectangle (3*\L, 2*\L);
				
	% roads with width > 0
	\draw[fill=black!20!white, draw=black!20!white]  (0, -\w/2) rectangle (\L, \w/2);
	\draw[fill=black!20!white, draw=black!20!white] 
		(\L, - \w/2) rectangle (2*\L, \w/2);
	\draw[fill=black!20!white, draw=black!20!white] 
		(\L - \w/2, 0) rectangle (\L + \w / 2, \L);
	\draw[fill=black!20!white, draw=black!20!white] 
		(2*\L - \w/2, 0) rectangle (2 * \L + \w / 2, \L);
	\draw[fill=black!20!white, draw=black!20!white] 
		(\L, \L - \w/2) rectangle (2 * \L, \L + \w/2);
	\draw[fill=black!20!white, draw=black!20!white] 
		(2*\L, \L - \w/2) rectangle (3 * \L, \L + \w/2);
		
	% overlapping parts of the roads
	\fill[black!20!white] (\L - \w / 2, 0) rectangle (\L + \w/2, \w/2);
	\fill[black!20!white] (2*\L - \w / 2, 0) rectangle (2*\L, \w/2);
	\fill[black!20!white] (2*\L - \w / 2, \L - \w/2) rectangle (2*\L + \w/2, \L);
	\fill[black!20!white] (\L, \L - \w /2) rectangle (\L + \w/2, \L);
		
	\draw[black!70!white] (0,-\L) grid[step=\h] (3*\L, 2*\L);
			
\end{tikzpicture}
}
		\subcaption{Schematic discretization.}
		\label{fig: discretization}
	\end{subfigure}
	\caption{The road network of the proof-of-concept example. The schematic spatial discretization includes the width of the roads.}
	\label{fig: proof-of-concept-example data}
\end{figure}
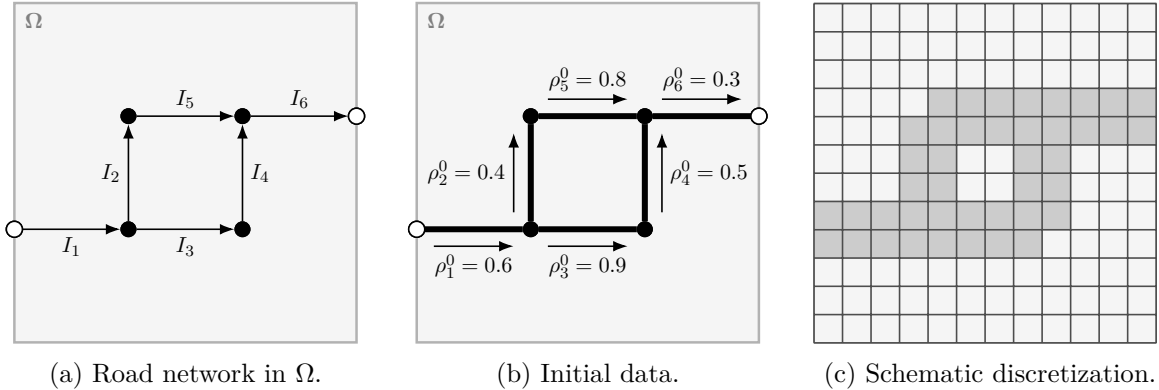

In this section, we present an overview of our model and discretization parameters (\cref{table: parameters}). The road network along with the initial data is shown in \cref{fig: road network,fig: initial data lwr} and the spatial discretization of the traffic model is illustrated in \cref{fig: discretization}.

In the following, we also comment briefly on the choice of the essential parameters.
\begin{table}[ht]
	\centering 
	\renewcommand{\arraystretch}{1.15}
	\caption{Overview of the parameter values defining the proof-of-concept example.}
	\label{table: parameters}
	\begin{tabular}[t]{ll}
		\toprule
		Parameter & Value\\
		\midrule
		\multicolumn{2}{c}{\textit{Traffic model}} \\
		\(T \) & 5 \\
		\(I_e\) & \(\left[0,1\right]\)\\		
		\(\rho_e^{\max}\) & 1  \\
		\(q_e^\inn\) & \( 0.25 \) \\
		\( \rho_e^0 \) & cf. \cref{fig: initial data lwr} \\
		\( \sigma_e\) & cf. \cref{fig: road network} \\
		\( \alpha_{j,i}, \) & \( 0.5 \)  \\
		\( \beta_{j,i}\) & \( 0.5 \) \\
		\bottomrule
	\end{tabular}\hspace{15pt}
	\begin{tabular}[t]{ll}
		\toprule
		Parameter & Value \\
		\midrule
		\multicolumn{2}{c}{\textit{Emission model}} \\
		\( \Omega \) & \([0,3] \times [0,3]\) \\
	%	\( \sigma_e\) & cf. \cref{fig: road network} \\
		\( w_e \) & \( 0.1 \) \\
		\( \theta \) & \( 0.5 \)  \\
		\multicolumn{2}{c}{\textit{Dispersion model}} \\
		\( \vv\) & \( \left( 1, 1 \right)^\top \)  \\
		\( \kappa \) & \( 0\) \\
		\( \mu \) & \( 10^{-6}\) \\
		\( \phi_0\) & 0  \\
		\bottomrule
	\end{tabular}\hspace{15pt}
	\begin{tabular}[t]{ll}
		\toprule
		Parameter & Value \\
		\midrule
		\multicolumn{2}{c}{\textit{Optimization}} \\
		\( \underline{V_e}\) & \( 0.25 \)  \\
		\( \overline{V_e}\) & \( 2 \)  \\
		\multicolumn{2}{c}{\textit{Discretization}} \\
		\( \Delta s, h\) & \( 0.05\)\\
		\( N_s\) & \( 20\) \\
		%	\( \nu_\CFL\) & \( 0.9\) \\
		\( N_h \) & \( 60 \)  \\
		\( \Delta t\) & \( 0.0083 \)  \\
		\( N_t \) & \(601\)  \\
		\bottomrule
	\end{tabular}
\end{table}

For the traffic model, we opt for the common academic parameters. This means all roads have a unit length of \( L_e^\lwr = 1 \) and the same maximal capacity of \( \rho_e^{\max} = 1\).
For the initial density distribution, we assume that the initial densities are constant along each road. 
Motivated by the academic scenario where typical speed limits are \( V_e^{\max} = 1\), we impose box constraints that bound the minimal and maximal allowed speed limit on a road as follows:
\begin{equation*}
	\S = \bigg \{ \V^{\max} \in \R^6  \mathrel{\Big|} \frac{1}{4} \le V_e \le 2 
	\text{ for all } e \in E \bigg \}.
\end{equation*}

Next, we suppose that in the dispersion model, the advection dominates the spread of air pollutants.
This assumption is reasonable because for gases a diffusion coefficient of around \( 10^{-6} \si{\square\metre\per\second} \) is typical, which is significantly smaller than a wind velocity of around \( 1 \si{\metre\per\second} \) that indicates calm air conditions.

Further, we assume that the control area \( \Omega \) is initially free of air pollutants, which implies
\begin{equation*}
	\phi(\x,0) = \phi_0(\x) = 0 \quad \text{for all} \quad \x \in \Omega.
\end{equation*}
Notice that this does not depict the real world. However, the alternative expression~\labelcref{eq: alternative expression J_diff} of the objective \( \J_\diff \) reveals that the term involving \( \phi_0 \) is independent of the speed-limits. In other words, the choice of \( \phi_0\) only results in a constant shift in \( \J_\diff \) and thus has no effect on the minimizers of \( \J_\diff\) and \( \J_\poll\).
Therefore, this assumption is unrestrictive.

To close the discussion on parameter selection, it is worth noting that the temporal step size \( \Delta t\) and the spatial step size \( h\) have to satisfy a CFL condition for the time-reversed adjoint equation, cf. \cref{remark cfl}, to ensures the stability of the numerical scheme.

After we have established the parameter choices, we analyze the solution of the optimization problem~\labelcref{eq: optimization problem} for different values of \( \delta\) next. Afterwards, we compare these optimal solutions to the one of \labelcref{eq: optimization problem 3D}.

\subsection{Results of the numerical experiments}
Before we continue with a detailed discussion of the results, let us comment on the choices of \( \delta\) which are used to showcase and compare optimal solutions. 

The parameter \( \delta \) regulates the influence of idle traffic on the caused air pollution. The choice \(\delta= 0 \) represents the case where idle traffic has no effect on the environment, while for \(\delta > 0\) idle traffic also contributes to air pollution, cf. \cref{eq: J_poll discrete}.
Hence, we consider \( \delta\in \{ 0, \frac{1}{2} \} \) to showcase the influence of idle traffic on air pollution and optimal speed limits.

\medbreak
We begin by solving the proof-of-concept example with \(\delta = 0 \). In this scenario the presence of a queue with idle vehicles has no impact on the environment. 

\begin{figure}[ht]
	\centering
	\begin{subfigure}[t]{.48\textwidth}
		\centering
\begin{tikzpicture}
	
	\begin{axis}[%
		width=0.72\textwidth, % 1
		height=0.567\textwidth, % 0.788
		at={(0\textwidth,0\textwidth)},
		scale only axis,
		xmin=0.3,
		xmax=1.025,
		xlabel style={font=\color{white!15!black}, font=\small},
		xlabel={$\J_\flow^h / \max \J^h_\flow$},
		ymin=0.975,
		ymax=1.7,
		ylabel style={font=\color{white!15!black}, font=\small},
		ylabel={$\J_\poll^h(\cdot;0) / \min \J_\poll^h(\cdot;0)$},
		ticklabel style = {font=\scriptsize},
		axis background/.style={fill=white},
		xmajorgrids,
		ymajorgrids,
		legend style={legend cell align=left, align=left, anchor=north west, 
			legend pos = north west, draw=white!15!black, font=\small}
		]
		\addplot[very thick, samples=50, smooth, domain=0:6, black, densely dashed, forget plot] 
		coordinates {(0.3,1)(1,1)};
		\addplot[very thick, samples=50, smooth, domain=0:6, black, densely dashed, forget plot] 
		coordinates {(1,1)(1, 3)};
		
		\addplot[only marks, mark=*, mark options={}, mark size=2pt, color=rwth, fill=rwth] table[row sep=crcr]{%
			x	y\\
			0.996481537621619	1.60007083772307\\
			0.512735152813021	1.0483788144392\\
			0.804158332562536	1.34297488317134\\
			0.792709814630393	1.32949148098035\\
			0.981990654023482	1.58866236792562\\
			0.40448638411715	1.00982088145845\\
			0.978302080500612	1.58758767128883\\
			0.66622674129452	1.13531195126905\\
			0.975586648200231	1.58585117039447\\
			0.629683011745276	1.11689960121747\\
			0.757058843529235	1.26958601199612\\
			0.693912239264176	1.18438184430408\\
			0.71566652113261	1.20552259327055\\
			0.37293431393404	1\\
			0.99944102379107	1.60299925111965\\
			0.443088061058373	1.02621192061857\\
			0.828266883105694	1.38278862141537\\
			0.904726530697175	1.50232208376818\\
			0.379456416758737	1.00586583514423\\
			0.895042764570956	1.48668153043681\\
			0.856704045216725	1.42676961097755\\
			0.763969226648526	1.27297448800451\\
			0.872590455250138	1.45285402144468\\
			0.734847864753792	1.23960790916043\\
			0.575694058939738	1.08999664795463\\
			0.473952793744308	1.03980419403566\\
			0.588215481122727	1.09625279169335\\
			0.615995991690132	1.11502922260605\\
			0.596534563701835	1.10177516441955\\
			0.696844596149391	1.18453134777158\\
			0.750879096073259	1.2605080696566\\
			0.516726904633838	1.04997757845211\\
			0.81120975699223	1.36213526487293\\
			0.865617205800092	1.4509972568448\\
			0.835403022488962	1.39012507882491\\
			0.453172489068096	1.03018723609763\\
			0.530833663206327	1.07603015634091\\
			0.86216137008243	1.44932881749533\\
			0.97377281744529	1.58540459138869\\
			0.489035126026451	1.04095704264708\\
			0.662005365099713	1.12504640597879\\
			0.957777223749107	1.56325569403817\\
			0.853652517467134	1.4246586080436\\
			0.545431403026903	1.07800710679062\\
			0.437039197916508	1.02490239248148\\
			0.959592026031151	1.56798212303917\\
			0.910789900371679	1.51310218740982\\
			0.928254057015316	1.53504362954491\\
			0.831855583909992	1.38931332060875\\
			0.918945290274139	1.52952019673836\\
			0.931198587846368	1.53684852916526\\
			0.591276426835301	1.09824406302359\\
			0.962996936972536	1.56816831827672\\
			0.880230592027725	1.4669981145884\\
			0.961260055374937	1.5681275819815\\
			0.964878991123463	1.5709109538264\\
			0.971985873890735	1.57966610092885\\
			0.942737022236928	1.5490623769004\\
			0.815764259058939	1.36400330247259\\
			0.966858121446068	1.57726918021509\\
			0.433416302987713	1.02416013410234\\
			0.526131051632785	1.05855653185348\\
			0.669045420443546	1.1363331656354\\
			0.380710835637856	1.00885596276035\\
			0.8200819308274	1.37029132053624\\
			0.417868188157339	1.0146430687961\\
			0.874260861698119	1.45512669913707\\
			0.968038962772814	1.5791055174209\\
			0.950659247441015	1.55996266189771\\
			0.985990126291584	1.59009055642323\\
			0.700545722955521	1.18515302377352\\
			0.947118175428474	1.55589364653819\\
			0.890505833092355	1.48241285821057\\
			0.94092356091853	1.54753954070014\\
			0.943465327803073	1.54948401591807\\
			0.952294569722368	1.56017485568427\\
			0.993050935856631	1.59290852199459\\
			0.956222996923705	1.56188708016552\\
			0.958474896449149	1.5636780948005\\
			0.944037588715441	1.55017484914849\\
		};
		\addlegendentry{Pareto Front}
		
		\addplot[only marks, mark=*, mark options={}, mark size=2pt, color=black, fill=black] table[row sep=crcr]{%
			x	y\\
			1	1\\
		};
		\addlegendentry{Ideal Vector}
		
	\end{axis}
\end{tikzpicture}%
	\subcaption{Pareto front and ideal vector.}
	\label{fig: pareto front without queue}
	\end{subfigure} \hfill
	\begin{subfigure}[t]{.48\textwidth}
		\centering
\begin{tikzpicture}
	
	\begin{axis}[%
		width=0.72\textwidth, % 1
		height=0.567\textwidth, % 0.788
		at={(0\textwidth,0\textwidth)},
		scale only axis,
		xmin=0.75,
		xmax=6.25,
		xtick={1,2,3,4,5,6},
		xticklabels={$I_1$, $I_2$, $I_3$, $I_4$, $I_5$, $I_6$},
		xlabel style={font=\color{white!15!black}, font=\small},
		xlabel={road},
		ymin=0.25,
		ymax=2.,
		ytick={0.25, 0.5, 0.75, 1, 1.25, 1.5, 1.75, 2},
		ylabel style={font=\color{white!15!black}, font=\small},
		ylabel={speed limit},
		ticklabel style = {font=\scriptsize},
		axis background/.style={fill=white},
		xmajorgrids,
		ymajorgrids
		]
		\addplot[only marks, mark=*, mark options={}, mark size=2pt, color=rwth, fill=rwth, forget plot] table[row sep=crcr]{%
			x	y\\
			1	0.25\\
			1	0.26171875\\
			1	0.37109375\\
			1	0.4140625\\
			1	0.484375\\
			1	0.5234375\\
			1	0.666015625\\
			1	0.6875\\
			1	0.734375\\
			1	0.76953125\\
			1	0.78125\\
			1	0.859375\\
			1	0.916015625\\
			1	0.947265625\\
			1	0.978515625\\
			1	0.984375\\
			1	1.015625\\
			1	1.041015625\\
			1	1.109375\\
			1	1.166015625\\
			1	1.234375\\
			1	1.26171875\\
			1	1.291015625\\
			1	1.359375\\
			1	1.484375\\
			1	1.609375\\
			1	1.734375\\
			1	2\\
		};
		\addplot[only marks, mark=*, mark options={}, mark size=2pt, color=rwth, fill=rwth, forget plot] table[row sep=crcr]{%
			x	y\\
			2	0.25\\
			2	0.796875\\
			2	1.1796875\\
			2	1.26171875\\
			2	1.4375\\
			2	1.47265625\\
			2	1.53515625\\
			2	1.72265625\\
			2	1.75390625\\
			2	1.810546875\\
			2	1.841796875\\
			2	1.84765625\\
			2	1.86328125\\
			2	1.904296875\\
			2	1.97265625\\
			2	2\\
		};
		\addplot[only marks, mark=*, mark options={}, mark size=2pt, color=rwth, fill=rwth, forget plot] table[row sep=crcr]{%
			x	y\\
			3	0.25\\
			3	0.34375\\
			3	0.46875\\
			3	0.59375\\
			3	1.09375\\
			3	1.103515625\\
			3	1.171875\\
			3	1.21484375\\
			3	1.421875\\
			3	1.46875\\
			3	1.478515625\\
			3	1.509765625\\
			3	1.541015625\\
			3	1.546875\\
			3	1.603515625\\
			3	1.6328125\\
			3	1.64453125\\
			3	1.671875\\
			3	1.69921875\\
			3	1.71484375\\
			3	1.728515625\\
			3	1.7421875\\
			3	1.890625\\
			3	1.90625\\
			3	2\\
		};
		\addplot[only marks, mark=*, mark options={}, mark size=2pt, color=rwth, fill=rwth, forget plot] table[row sep=crcr]{%
			x	y\\
			4	0.34375\\
			4	0.46875\\
			4	0.5625\\
			4	0.6328125\\
			4	0.6875\\
			4	0.7421875\\
			4	1.015625\\
			4	1.0625\\
			4	1.09765625\\
			4	1.1328125\\
			4	1.25\\
			4	1.2890625\\
			4	1.3125\\
			4	1.34375\\
			4	1.400390625\\
			4	1.431640625\\
			4	1.4375\\
			4	1.46875\\
			4	1.4765625\\
			4	1.494140625\\
			4	1.5625\\
			4	1.58984375\\
			4	1.619140625\\
			4	1.6875\\
			4	1.7265625\\
			4	1.8515625\\
			4	1.86328125\\
			4	1.994140625\\
			4	2\\
		};
		\addplot[only marks, mark=*, mark options={}, mark size=2pt, color=rwth, fill=rwth, forget plot] table[row sep=crcr]{%
			x	y\\
			5	0.25\\
			5	0.26171875\\
			5	0.359375\\
			5	0.40625\\
			5	0.5\\
			5	0.734375\\
			5	0.75\\
			5	0.953125\\
			5	1.015625\\
			5	1.234375\\
			5	1.25\\
			5	1.34375\\
			5	1.453125\\
			5	1.603515625\\
			5	1.671875\\
			5	1.69921875\\
			5	1.76171875\\
			5	1.796875\\
			5	1.80859375\\
			5	1.82421875\\
			5	1.84375\\
			5	1.853515625\\
			5	2\\
		};
		\addplot[only marks, mark=*, mark options={}, mark size=2pt, color=rwth, fill=rwth, forget plot] table[row sep=crcr]{%
			x	y\\
			6	2\\
		};
	\end{axis}
\end{tikzpicture}%
		\subcaption{Values of the Pareto-optimal speed limit policies.}
		\label{fig: speed limit range without queue}
	\end{subfigure}
	\caption{The optimal solution of the proof-of-concept example with \(\delta = 0\), i.e., 
		air pollution is sorely estimated by active traffic on the road network.}
	\label{fig: problem without queue}
\end{figure}
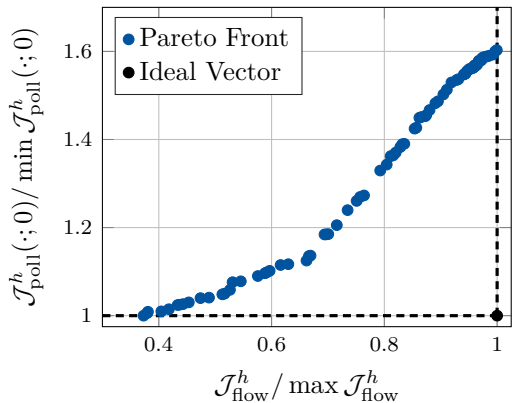
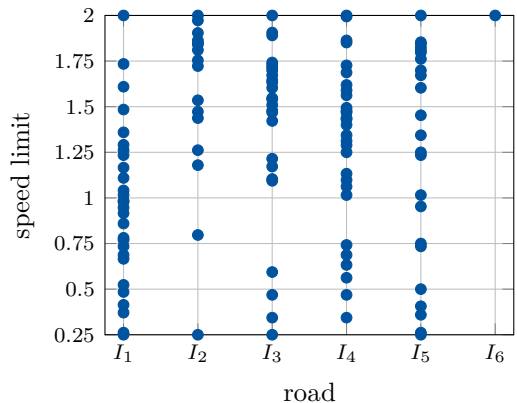

In \cref{fig: problem without queue}, we visualize the resulting optimal solution, including the Pareto front in \cref{fig: pareto front without queue}.
Notice that the ideal vector and the Pareto front are located in the lower right corner of the plot instead of the lower left corner because we maximize \( \J_\flow^h\).
Furthermore, the axes of \cref{fig: pareto front without queue} are normalized by the optimal values of the objectives, i.e., the accumulated traffic flow and the indicator to measure the caused air pollution are divided by
\begin{equation*}
	\max_{\V^{\max} \in \S} \J_\flow^h(\V^{\max}) \quad \text{and} \quad
	\min_{\V^{\max} \in \S} \J_\poll^h(\V^{\max};\delta),
\end{equation*}
respectively. Hence, the axes represent the relative difference of the function values compared to the optimal value, which is given by one as the ideal vector becomes \( (1,1)^\top \) for the normalized objectives.

For example, \cref{fig: pareto front without queue} shows that sorely minimizing the contribution to air pollution results in the efficient point \( (0.37,1)^\top \). This point translates to: while the contribution to air pollution is minimal, the accumulated traffic flow decreases by \( 63 \% \) compared to its ideal value.
The reversed case yields the efficient point \( (1,1.6)^\top \) where we maximize the economic efficiency. Here the environmental impact increases by \( 60 \% \) compared to its ideal value while the accumulated traffic flow attains its maximum.

An example for an efficient point where neither of the objectives attains its optimum, is the point \( (0.8, 1.34)^\top \). Here the environmental impact is increased by \( 34 \% \) while the economic efficiency decreases by \( 20 \% \).

\medbreak
Because each efficient solution of the Pareto front corresponds to a Pareto-optimal speed limit policy, we investigate the Pareto-optimal speed limit policies in more detail now.

\cref{fig: speed limit range without queue} shows the computed Pareto-optimal speed limits of the individual roads.
Notice this provides information on the range of the optimal speed limits of a road but not on the concrete optimal policies since those are vectors, i.e., randomly connecting the optimal speed limits of the individual roads illustrated in \cref{fig: speed limit range without queue}, does not yield an optimal policy.
Nevertheless, \cref{fig: speed limit range without queue} reveals that the two objectives are conflicting as the range of the optimal speed limits nearly covers the entire interval \( [0.25,2] \) for each road, except for the exit road \( e = 6 \).
At the exit road, the two objectives \textquote{agree} on an optimal speed limit, i.e., \( V_6^{\max} = 2\). 
A possible explanation may be that, as we impose free flow at the exit road and the initial density of this  road is relatively low, the vehicles should and can leave the network as fast as possible because then they do not contribute to the air pollution anymore -- at least from the modeling point of view.
Further, the flux function scales linearly with the speed limit, i.e., for a fixed density, the traffic flow increases as the speed limit increases, resulting in a higher traffic flow.

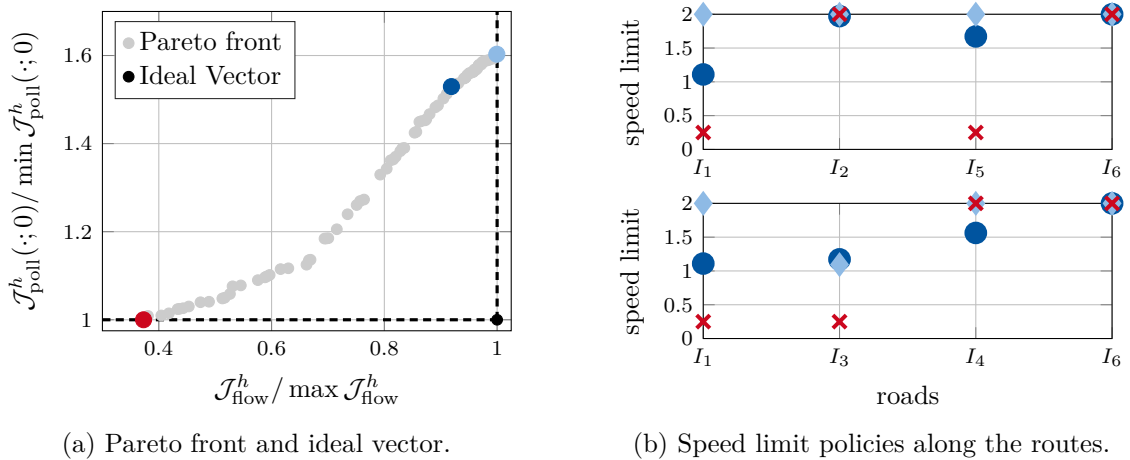
\begin{figure}[ht]
	\centering
	\begin{subfigure}[t]{.48\textwidth}
		\centering
\begin{tikzpicture}
	
	\begin{axis}[%
		width=0.72\textwidth, % 1
		height=0.567\textwidth, % 0.788
		at={(0\textwidth,0\textwidth)},
		scale only axis,
		xmin=0.3,
		xmax=1.025,
		xlabel style={font=\color{white!15!black}, font=\small},
		xlabel={$\J_\flow^h / \max \J^h_\flow$},
		ymin=0.975,
		ymax=1.7,
		ylabel style={font=\color{white!15!black}, font=\small},
		ylabel={$\J_\poll^h(\cdot;0) / \min \J_\poll^h(\cdot;0)$},
		ticklabel style = {font=\scriptsize},
		axis background/.style={fill=white},
		xmajorgrids,
		ymajorgrids,
		legend style={legend cell align=left, align=left, anchor=north west, 
			legend pos = north west, draw=white!15!black, font=\small}
		]
		\addplot[very thick, samples=50, smooth, domain=0:6, black, densely dashed, forget plot] 
		coordinates {(0.3,1)(1,1)};
		\addplot[very thick, samples=50, smooth, domain=0:6, black, densely dashed, forget plot] 
		coordinates {(1,1)(1, 3)};
		
		\addplot[only marks, mark=*, mark options={}, mark size=2pt, color=black!20!white, fill=black!20!white] table[row sep=crcr]{%
			x	y\\
			0.996481537621619	1.60007083772307\\
			0.512735152813021	1.0483788144392\\
			0.804158332562536	1.34297488317134\\
			0.792709814630393	1.32949148098035\\
			0.981990654023482	1.58866236792562\\
			0.40448638411715	1.00982088145845\\
			0.978302080500612	1.58758767128883\\
			0.66622674129452	1.13531195126905\\
			0.975586648200231	1.58585117039447\\
			0.629683011745276	1.11689960121747\\
			0.757058843529235	1.26958601199612\\
			0.693912239264176	1.18438184430408\\
			0.71566652113261	1.20552259327055\\
			0.37293431393404	1\\
			0.99944102379107	1.60299925111965\\
			0.443088061058373	1.02621192061857\\
			0.828266883105694	1.38278862141537\\
			0.904726530697175	1.50232208376818\\
			0.379456416758737	1.00586583514423\\
			0.895042764570956	1.48668153043681\\
			0.856704045216725	1.42676961097755\\
			0.763969226648526	1.27297448800451\\
			0.872590455250138	1.45285402144468\\
			0.734847864753792	1.23960790916043\\
			0.575694058939738	1.08999664795463\\
			0.473952793744308	1.03980419403566\\
			0.588215481122727	1.09625279169335\\
			0.615995991690132	1.11502922260605\\
			0.596534563701835	1.10177516441955\\
			0.696844596149391	1.18453134777158\\
			0.750879096073259	1.2605080696566\\
			0.516726904633838	1.04997757845211\\
			0.81120975699223	1.36213526487293\\
			0.865617205800092	1.4509972568448\\
			0.835403022488962	1.39012507882491\\
			0.453172489068096	1.03018723609763\\
			0.530833663206327	1.07603015634091\\
			0.86216137008243	1.44932881749533\\
			0.97377281744529	1.58540459138869\\
			0.489035126026451	1.04095704264708\\
			0.662005365099713	1.12504640597879\\
			0.957777223749107	1.56325569403817\\
			0.853652517467134	1.4246586080436\\
			0.545431403026903	1.07800710679062\\
			0.437039197916508	1.02490239248148\\
			0.959592026031151	1.56798212303917\\
			0.910789900371679	1.51310218740982\\
			0.928254057015316	1.53504362954491\\
			0.831855583909992	1.38931332060875\\
			0.918945290274139	1.52952019673836\\
			0.931198587846368	1.53684852916526\\
			0.591276426835301	1.09824406302359\\
			0.962996936972536	1.56816831827672\\
			0.880230592027725	1.4669981145884\\
			0.961260055374937	1.5681275819815\\
			0.964878991123463	1.5709109538264\\
			0.971985873890735	1.57966610092885\\
			0.942737022236928	1.5490623769004\\
			0.815764259058939	1.36400330247259\\
			0.966858121446068	1.57726918021509\\
			0.433416302987713	1.02416013410234\\
			0.526131051632785	1.05855653185348\\
			0.669045420443546	1.1363331656354\\
			0.380710835637856	1.00885596276035\\
			0.8200819308274	1.37029132053624\\
			0.417868188157339	1.0146430687961\\
			0.874260861698119	1.45512669913707\\
			0.968038962772814	1.5791055174209\\
			0.950659247441015	1.55996266189771\\
			0.985990126291584	1.59009055642323\\
			0.700545722955521	1.18515302377352\\
			0.947118175428474	1.55589364653819\\
			0.890505833092355	1.48241285821057\\
			0.94092356091853	1.54753954070014\\
			0.943465327803073	1.54948401591807\\
			0.952294569722368	1.56017485568427\\
			0.993050935856631	1.59290852199459\\
			0.956222996923705	1.56188708016552\\
			0.958474896449149	1.5636780948005\\
			0.944037588715441	1.55017484914849\\
		};
		\addlegendentry{Pareto front}
		
		\addplot[only marks, mark=*, mark options={}, mark size=2pt, color=black, fill=black] table[row sep=crcr]{%
			x	y\\
			1	1\\
		};
		\addlegendentry{Ideal Vector}
		
		\addplot[only marks, mark=*, mark options={}, mark size=3pt, color=rwth, fill=rwth, forget plot] table[row sep=crcr]{%
			x	y\\
			0.918945290274139	1.52952019673836\\
		};
		%\addlegendentry{middle}
		
		\addplot[only marks, mark=*, mark options={}, mark size=3pt, color=rwth-50, fill=rwth-50, forget plot] table[row sep=crcr]{%
			x	y\\
			0.99944102379107	1.60299925111965\\
		};
		%\addlegendentry{flow}
		
		\addplot[only marks, mark=*, mark options={}, mark size=3pt, color=rot, fill=rot, forget plot] table[row sep=crcr]{%
			x	y\\
			0.37293431393404	1\\
		};
		%\addlegendentry{eco}
		
	\end{axis}
\end{tikzpicture}%
	\subcaption{Pareto front and ideal vector.}
	\label{fig: pareto front highlights}
	\end{subfigure} \hfill
	\begin{subfigure}[t]{.48\textwidth}
		\centering
\begin{tikzpicture}
	% upper route
	\begin{axis}[%
		name=upperRoute,
		width=0.72\textwidth, %0.45
		height=0.24\textwidth, % 0.15
		at={(0\textwidth,0\textwidth)},
		scale only axis,
		xmin=1,
		xmax=4,
		xtick={1,2,3,4},
		xticklabels={$I_1$, $I_2$, $I_5$, $I_6$},
		xlabel style={font=\color{white!15!black}, font=\small},
		ymin=0,
		ymax=2,
		ylabel style={font=\color{white!15!black}, font=\small},
		ylabel={speed limit},
		axis background/.style={fill=white},
		ticklabel style = {font=\scriptsize},
		xmajorgrids,
		ymajorgrids,
		]
		\addplot[only marks, line width=1.5pt, color=rwth, mark=*, mark options={solid, rwth}, mark size = 3.5pt, forget plot]
		table[row sep=crcr]{%
			1	1.109375\\
			2	1.97265625\\
			3	1.671875\\
			4	2\\
		};
		%\addlegendentry{middle}
				
		\addplot[only marks, line width=1.5pt,color=rwth-50, mark=diamond*, mark options={solid, rwth-50}, mark size = 3.5pt, forget plot]
		table[row sep=crcr]{%
			1	2\\
			2	2\\
			3	2\\
			4	2\\
		};
		%\addlegendentry{max flow}
		
			\addplot[only marks, line width=1.5pt, color=rot, mark=x, mark options={solid, rot}, mark size = 3.5pt, forget plot]
		table[row sep=crcr]{%
			1	0.25\\
			2	2\\
			3	0.25\\
			4	2\\
		};
		%\addlegendentry{min eco}

	\end{axis}
	
	% lower route!!
	\begin{axis}[%
		width=0.72\textwidth, %0.45
		height=0.24\textwidth, % 0.15
		at=(upperRoute.below south west), anchor=above north west,
		scale only axis,
		xmin=1,
		xmax=4,
		xtick={1,2,3,4},
		xticklabels={$I_1$, $I_3$, $I_4$, $I_6$},
		ymin=0.2,
		ymax=2,
		xlabel style={font=\color{white!15!black}, font=\small},
		xlabel={roads},
		ymin=0,
		ymax=2,
		ylabel style={font=\color{white!15!black}, font=\small},
		ylabel={speed limit},
		axis background/.style={fill=white},
		ticklabel style = {font=\scriptsize},
		xmajorgrids,
		ymajorgrids,
		]
			\addplot [only marks, color=rwth, mark=*, mark options={solid, rwth},
		line width=1.5pt, mark size = 3.5pt, forget plot]
		table[row sep=crcr]{%
			1	1.109375\\
			2	1.171875\\
			3	1.5625\\
			4	2\\
		};
		%\addlegendentry{middle}
		\addplot  [only marks, color=rwth-50, mark=diamond*, mark options={solid, rwth-50}, 
		line width=1.5pt, mark size = 3.5pt, forget plot]
		table[row sep=crcr]{%
			1	2\\
			2	1.09375\\
			3	2\\
			4	2\\
		};
		%\addlegendentry{max flow}
		
		\addplot[only marks, color=rot, mark=x, mark options={solid, rot}, 
		line width=1.5pt, mark size = 3.5pt, forget plot]
		table[row sep=crcr]{%
			1	0.25\\
			2	0.25\\
			3	2\\
			4	2\\
		};
		%\addlegendentry{min eco}
		
	\end{axis}
\end{tikzpicture}%
	\subcaption{Speed limit policies along the routes.}
	\label{fig: speed limit routes without queue}
	\end{subfigure}
	\caption{Pareto-optimal speed limit policies along the two routes of the sample road network.
	The correspondence between efficient solutions and the Pareto-optimal speed limit policies is
	color-coded (in red, dark blue and light blue).}
	\label{fig: routes without queue}
\end{figure}

%\medbreak
On the other hand, we visualize certain speed limit policies along the two routes of the road network in \cref{fig: routes without queue} to gain an impression of what optimal speed limit policies look like and if they impose any difficulties considering real-world applications. 
The two routes a driver can take are the \textquote{lower} and \textquote{upper} route, which are given by
\begin{equation*}
	I_1 \rightarrow I_3 \rightarrow I_4 \rightarrow I_6 
	\quad \text{and} \quad
	I_1 \rightarrow I_2 \rightarrow I_5 \rightarrow I_6,
\end{equation*}
respectively, see \cref{fig: road network}.

The speed limits policy corresponding to the efficient solution, where \( \J_\poll^h \) attains its minimum, is a bang-bang control, i.e., the speed limits are either given by the lower or upper bound of the inequality constraint.
Along the upper route, the speed limits alternate between the lower and upper bound.
Even though this is unproblematic from the mathematical point of view, it may be in real-world applications because in reality, drivers may violate the lower speed limit when coming from a high-speed-limit road.

A similar problem occurs at the diverging junction, where the speed limits of the outgoing roads are \( V_2^{\max} = 0.25 \) and \( V_3^{\max} = 2 \). Again, this does not impose a mathematical problem, but in the real world, drivers can \textquote{choose} between the outgoing roads, and hence, a majority may opt for road \( e = 3 \) because it allows higher speeds. 
However, at each diverging junction, we prescribe speed-limit-independent distribution rates to describe the percentage of incoming drivers accessing an outgoing road.
Hence, a higher speed limit creates an incentive to opt for the corresponding road, which then results in a violation of the modeled distribution rates, i.e., the model does not depict the reality anymore.
Both issues present a subject of future research: adjusting the constraints on the speed limit policies to avoid a big difference between consecutive speed limits or speed limits of outgoing roads at a junction.

Another speed limit policy in \cref{fig: routes without queue} is marked in dark blue, where neither of the objectives attains its optimum.
This speed limit policy shows a similar behavior compared to the one corresponding to minimal environmental impact, but the increase of the speed limits along the lower route is shallower. 
Further, the lowest speed limit along the upper route is around \( V_e^{\max} \approx 0.6\) compared to \( V_e^{\max} = 0.25 \) for the example of minimal environmental impact.

\begin{figure}[ht]
	\centering
	\begin{subfigure}[t]{.48\textwidth}
		\centering
\begin{tikzpicture}
	
	\begin{axis}[%
		width=0.72\textwidth, % 1
		height=0.567\textwidth, % 0.788
		at={(0\textwidth,0\textwidth)},
		scale only axis,
		xmin=0.5,
		xmax=1.025,
		xlabel style={font=\color{white!15!black}, font=\small},
		xlabel={$\J_\flow^h / \max \J^h_\flow$},
		ymin=0.985,
		ymax=1.25,
		ylabel style={font=\color{white!15!black}, font=\small},
		ylabel={$\J_\poll^h\left(\cdot;\frac{1}{2}\right) / \min \J_\poll^h\left(\cdot;\frac{1}{2}\right)$},
		ticklabel style = {font=\scriptsize},
		axis background/.style={fill=white},
		xmajorgrids,
		ymajorgrids,
		legend style={legend cell align=left, align=left, anchor=north west, 
			legend pos = north west, draw=white!15!black, font=\small}
		]
		\addplot[very thick, samples=50, smooth, domain=0:6, black, densely dashed, forget plot] 
		coordinates {(0.5,1)(1,1)};
		\addplot[very thick, samples=50, smooth, domain=0:6, black, densely dashed, forget plot] 
		coordinates {(1,1)(1, 1.25)};
		
		\addplot[only marks, mark=*, mark options={}, mark size=2pt, color=rwth, fill=rwth] table[row sep=crcr]{%
			x	y\\
		0.99510273834397	1.21114920692553\\
		0.53762642490493	1.02321267961985\\
		0.775561969283996	1.08871852948844\\
		0.573953367743057	1.03183872868223\\
		0.614001429010663	1.03848859296699\\
		0.510325305496653	1.00781716569255\\
		0.825217398217446	1.11435936315716\\
		0.656580902646739	1.04547724988731\\
		0.75455923410347	1.07992430490492\\
		0.602368765607133	1.03706027827905\\
		0.910925149833243	1.16164139434448\\
		0.705012052424752	1.06184939147516\\
		0.516377112874728	1.01498446359441\\
		0.719769827394763	1.06858678946639\\
		0.501705400930418	1.00428822416443\\
		0.63254290539619	1.03918927102754\\
		0.764775338340774	1.08405002926337\\
		0.802583117508615	1.10092917891779\\
		0.55793274980429	1.02908145612441\\
		0.866992466923337	1.13737617869665\\
		0.784329997368422	1.093347438827\\
		0.723965285297225	1.06937489943054\\
		0.693202656357988	1.05738156246614\\
		0.957777223749107	1.18998790738076\\
		0.740809072982964	1.0737133604106\\
		0.959592026031151	1.19336921946218\\
		0.587023905041016	1.03417948388228\\
		0.967132257892364	1.19625200702155\\
		0.810335403853001	1.10664120551441\\
		0.929371793053975	1.16967268570874\\
		0.797682196303834	1.09958039760601\\
		0.876481232550363	1.14612041621154\\
		0.785631489199676	1.09769577665659\\
		0.936475920694774	1.17470223761731\\
		0.703549265850537	1.06138427044853\\
		0.970205776268762	1.20223703934078\\
		0.681400760403829	1.05316344228687\\
		0.733488515636229	1.07199650546533\\
		0.718697717971542	1.06780112025372\\
		0.923217611906895	1.16705369436106\\
		0.907782160244926	1.15554925226457\\
		0.643772152179466	1.04176961271101\\
		0.962996936972536	1.19353824344369\\
		0.940069804520093	1.17697479586527\\
		0.826431755760132	1.11492318310376\\
		0.942737022236929	1.17976692578507\\
		0.660583948822576	1.04703062165068\\
		0.814974730966834	1.10901857070312\\
		0.848642104919183	1.12812190854065\\
		0.658017302040712	1.04579495564323\\
		0.918945290274139	1.1656235708839\\
		0.805524528728601	1.10356279757389\\
		0.610030470283887	1.037213359993\\
		0.983137504495652	1.20796738525691\\
		0.77823849989148	1.0892017408783\\
		0.954047714942694	1.18983898523607\\
		0.828493227959402	1.11523206253734\\
		0.830368798840913	1.11609917954556\\
		0.97631660607409	1.20548955086793\\
		0.950659247441015	1.18746116317617\\
		0.715510863385365	1.0636301466468\\
		0.71197711468477	1.06262647820891\\
		0.820316477012666	1.11301058184537\\
		0.899061692736782	1.15315423973947\\
		0.861390261433667	1.13393484053062\\
		0.646225047481246	1.04211244041263\\
		0.871031702205886	1.14175174152026\\
		0.903553273136339	1.15525427201236\\
		0.858137530187014	1.13336198904928\\
		0.883728918710059	1.14977612607193\\
		0.984106995448645	1.20895122988352\\
		0.669847361676084	1.04990339771631\\
		0.842269714403721	1.12430215721511\\
		0.910789900371679	1.15687775099549\\
		0.858995126211545	1.13385065878194\\
		0.854732004348756	1.13101183542998\\
		0.673394464549666	1.05070799057241\\
		0.947118175428475	1.18439500605834\\
		0.954235700455751	1.18989865405105\\
		0.833135431800232	1.11920350597193\\
		};
		\addlegendentry{Pareto Front}
		
		\addplot[only marks, mark=*, mark options={}, mark size=2pt, color=black, fill=black] table[row sep=crcr]{%
			x	y\\
			1	1\\
		};
		\addlegendentry{Ideal Vector}
		
	\end{axis}
\end{tikzpicture}%
	\subcaption{Pareto front and ideal vector}
	\label{fig: pareto front with queue}
	\end{subfigure} \hfill
	\begin{subfigure}[t]{.48\textwidth}
		\centering
\begin{tikzpicture}
	
	\begin{axis}[%
		width=0.72\textwidth, % 1
		height=0.567\textwidth, % 0.788
		at={(0\textwidth,0\textwidth)},
		scale only axis,
		xmin=0.75,
		xmax=6.25,
		xtick={1,2,3,4,5,6},
		xticklabels={$I_1$, $I_2$, $I_3$, $I_4$, $I_5$, $I_6$},
		xlabel style={font=\color{white!15!black}, font=\small},
		xlabel={road},
		ymin=0.2,
		ymax=2.05,
		ytick={0.25, 0.5, 0.75, 1, 1.25, 1.5, 1.75, 2},
		ylabel style={font=\color{white!15!black}, font=\small},
		ylabel={speed limit},
		ticklabel style = {font=\scriptsize},
		axis background/.style={fill=white},
		xmajorgrids,
		ymajorgrids
		]
		\addplot[only marks, mark=*, mark options={}, mark size=2pt, color=rwth, fill=rwth, forget plot] table[row sep=crcr]{%
			x	y\\
			1	1.041015625\\
			1	1.046875\\
			1	1.048828125\\
			1	1.09765625\\
			1	1.109375\\
			1	1.125\\
			1	1.15234375\\
			1	1.166015625\\
			1	1.19921875\\
			1	1.23046875\\
			1	1.234375\\
			1	1.26171875\\
			1	1.275390625\\
			1	1.359375\\
			1	1.37890625\\
			1	1.4140625\\
			1	1.5625\\
			1	1.666015625\\
			1	2\\
		};
		\addplot[only marks, mark=*, mark options={}, mark size=2pt, color=rwth, fill=rwth, forget plot] table[row sep=crcr]{%
			x	y\\
			2	0.568359375\\
			2	0.578125\\
			2	0.98828125\\
			2	1.015625\\
			2	1.04296875\\
			2	1.33984375\\
			2	1.58984375\\
			2	1.701171875\\
			2	1.70703125\\
			2	1.7265625\\
			2	1.779296875\\
			2	1.83984375\\
			2	1.84765625\\
			2	1.904296875\\
			2	1.90625\\
			2	1.91015625\\
			2	1.97265625\\
			2	2\\
		};
		\addplot[only marks, mark=*, mark options={}, mark size=2pt, color=rwth, fill=rwth, forget plot] table[row sep=crcr]{%
			x	y\\
			3	0.25\\
			3	0.359375\\
			3	0.400390625\\
			3	0.44140625\\
			3	0.48046875\\
			3	0.54296875\\
			3	0.796875\\
			3	0.85546875\\
			3	1.041015625\\
			3	1.046875\\
			3	1.109375\\
			3	1.171875\\
			3	1.330078125\\
			3	1.353515625\\
			3	1.41015625\\
			3	1.421875\\
			3	1.478515625\\
			3	1.53515625\\
			3	1.5625\\
			3	1.603515625\\
			3	1.671875\\
			3	1.69921875\\
			3	1.78125\\
			3	1.796875\\
			3	1.80859375\\
			3	2\\
		};
		\addplot[only marks, mark=*, mark options={}, mark size=2pt, color=rwth, fill=rwth, forget plot] table[row sep=crcr]{%
			x	y\\
			4	0.25\\
			4	1\\
			4	1.166015625\\
			4	1.220703125\\
			4	1.369140625\\
			4	1.4375\\
			4	1.494140625\\
			4	1.5\\
			4	1.5625\\
			4	1.58984375\\
			4	1.6875\\
			4	1.7265625\\
			4	1.75390625\\
			4	1.80859375\\
			4	1.86328125\\
			4	2\\
		};
		\addplot[only marks, mark=*, mark options={}, mark size=2pt, color=rwth, fill=rwth, forget plot] table[row sep=crcr]{%
			x	y\\
			5	0.25\\
			5	0.26171875\\
			5	0.27734375\\
			5	0.32421875\\
			5	0.35546875\\
			5	0.359375\\
			5	0.375\\
			5	0.38671875\\
			5	0.3984375\\
			5	0.5\\
			5	0.51171875\\
			5	0.57421875\\
			5	0.60546875\\
			5	0.69921875\\
			5	0.76171875\\
			5	0.82421875\\
			5	0.837890625\\
			5	0.88671875\\
			5	0.94921875\\
			5	1.01171875\\
			5	1.103515625\\
			5	1.228515625\\
			5	1.26171875\\
			5	1.3984375\\
			5	1.546875\\
			5	1.603515625\\
			5	1.609375\\
			5	1.671875\\
			5	1.69921875\\
			5	1.783203125\\
			5	1.796875\\
			5	1.86328125\\
			5	2\\
		};
		\addplot[only marks, mark=*, mark options={}, mark size=2pt, color=rwth, fill=rwth, forget plot] table[row sep=crcr]{%
			x	y\\
			6	2\\
		};
	\end{axis}
\end{tikzpicture}%
		\subcaption{Values of the Pareto-optimal speed limit policies}
		\label{fig: speed limit range with queue}
	\end{subfigure}
	\caption{The optimal solution of the proof-of-concept example with \( \delta = 1/2\), i.e., 
		air pollution is estimated by active traffic on the road network 
		and idle traffic due to the queues.}
	\label{fig: problem with queue}
\end{figure}
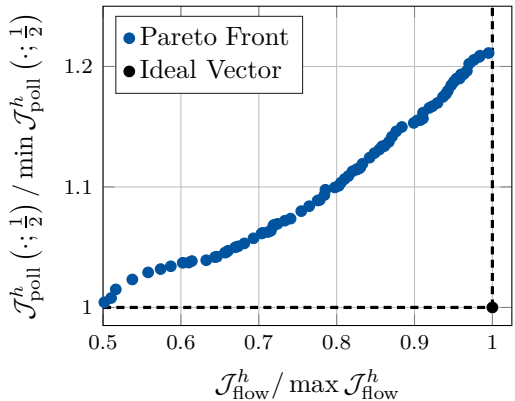
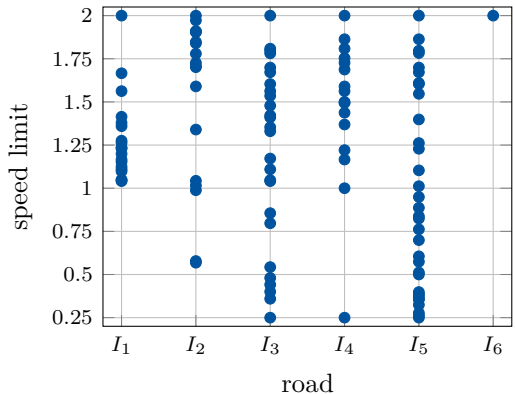

\medbreak
Next, we solve the proof-of-concept example with \(\delta = 1/2 \), i.e., the indicator measuring the environmental impact of vehicular traffic also accounts for idle traffic. 
\cref{fig: problem with queue} visualizes the resulting optimal solution consisting of the Pareto front (\cref{fig: pareto front with queue}) and the Pareto-optimal speed limit policies (\cref{fig: speed limit range with queue}).
Again, the axes of the Pareto front are normalized by the optimal values of the objectives.

The Pareto front reveals that for this modeling approach maximizing the accumulated traffic flow results in an increase of the environmental impact by around \( 20 \% \) while minimizing the environmental impact decreases the flow by around \( 50 \% \) compared to its optimal value.
An efficient point where neither of the objectives is ideal is, for example, the point \( (0.8,1.1)^\top \).

The Pareto-optimal speed limits policies are similar in terms of the range compared to the modeling approach where \( \delta= 0 \), i.e., they also indicate that the objectives are conflicting.
However, we observe a subtle change in the range of the speed limits of the access road: The lowest set speed limit is \( V_1^{\max} \approx 1 \) compared to \( V_1^{\max} = 0.25 \) for \( \delta = 0 \).
We can attribute this phenomenon to the idle traffic accounted by the indicator to measure air pollution caused by vehicular traffic.
The indicator accounts for the idle traffic by adding the average accumulated queue length to \( \J^h_\diff \), which then creates an incentive to reduce the queue length.
The correspondence between the queue length and the speed limits is the following:
The supply of a road scales linearly with its speed limit, i.e., an increased speed limit results in a higher flow that can access the road. The more flow can access the road network, the less excess accrues, which results in a shorter queue.

We omit a discussion on specific speed limits policies if \( \delta = 1/2 \) as this would reveal the same issues that can be observed in \cref{fig: speed limit routes without queue}, where the modeling neglects the contribution of the idle traffic to air pollution.

\begin{figure}[ht]
	\centering
	\begin{subfigure}[t]{.49\textwidth}
		\centering
		\input{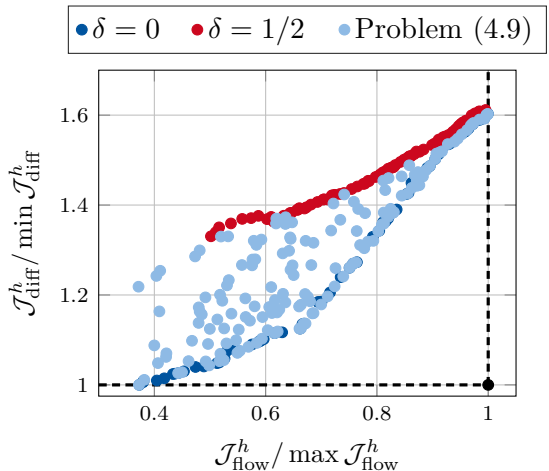}
		\subcaption{\( (\J^h_\flow, \J_\poll^h(\cdot;0)) \) coordinates.}
		\label{fig: coordinates without queue}
	\end{subfigure} \hfill
	\begin{subfigure}[t]{.49\textwidth}
		\centering
		\input{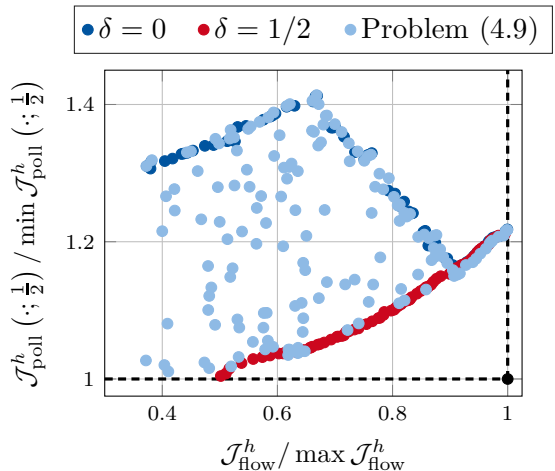}
		\subcaption{\( (\J^h_\flow, \J_\poll^h(\cdot;1/2)) \) coordinates.}
		\label{fig: coordinates with queue}
	\end{subfigure}
	\caption{Comparison of the Pareto fronts of the optimization problem~\labelcref{eq: optimization problem} for different values of \( \delta \), and of~\labelcref{eq: optimization problem 3D} by representing them in
		different coordinate systems.}
	\label{fig: coordinates}
\end{figure}

\begin{figure}[ht]
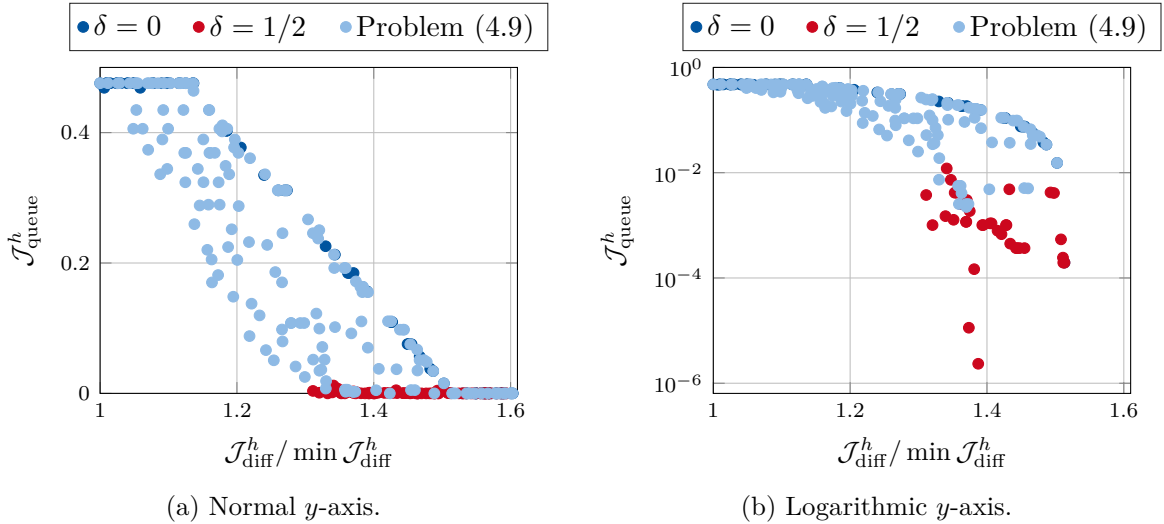

	\centering
	\begin{subfigure}[t]{.49\textwidth}
		\centering
		\input{./tikz_revision/composition_environmental_impact_review}
		\subcaption{Normal \(y\)-axis.}
	\end{subfigure} \hfill
	\begin{subfigure}[t]{.49\textwidth}
		\centering
		\input{./tikz_revision/composition_environmental_impact_loglog_review}
		\subcaption{Logarithmic \(y\)-axis.}
	\end{subfigure}
	\caption{Comparison of the contribution to air pollution caused by active traffic 
	(measured by \( \J^h_\diff \)) and idle traffic (measurement based on the average queue length). Both plots show the same image but the \(y\)-axis is scaled differently.}
	\label{fig: composition}
\end{figure}

\medbreak
At last, we compare the optimal solutions to the two modeling approaches~\labelcref{eq: optimization problem,eq: optimization problem 3D}.
Therefore, we display the computed Pareto fronts in \( (\J_\flow^h,\J_\poll^h(\cdot; \delta)) \) coordinates for \( \delta \in \{0, 1/2\} \), see \cref{fig: coordinates}, and in  \( (\J_\diff^h,\J_\queue^h) \) coordinates, see \cref{fig: composition}. Notice, that we do not normalize the axis \( \J_\queue^h\) as the minimal value of \( \J_\queue^h\) is in fact zero.

\Cref{fig: coordinates} shows that each Pareto front is Pareto-optimal in \textquote{its} coordinate system. If the accumulated traffic flow attains values close to its maximal value, the two fronts coincide for \(\delta= 0 \) and \( \delta = 1/2\).
Furthermore, we observe that the optimal solution to the optimization problem~\labelcref{eq: optimization problem 3D} involving three objectives, (a) lies between the two Pareto fronts involving only two, and (b) seems to have a tendency to share more efficient points with the solution to the problem~\labelcref{eq: optimization problem} with \(\delta = 0\).
In \cref{fig: coordinates without queue}, we observe that the optimal solution corresponding to \( \delta= 1/2 \) allows a higher traffic flow in general and that the active traffic contributes more to air pollution compared to the optimal solution corresponding to \( \delta = 0\).
A possible cause may be the penalty for the presence of idle traffic in the case \( \delta = 1/2 \), which then allows \textquote{more} active traffic on the network in the attempt to reduce the contribution of idle traffic.

Further, in \cref{fig: coordinates with queue}, we see that the combined contribution of active and idle traffic to air pollution is higher in case \(\delta= 0 \) compared to the case \( \delta = 1/2 \).
Again, a possible explanation may be that the indicator measuring the environmental impact also accounts for the idle traffic if \( \delta = 1/2 \), while the case  \( \delta = 0 \) neglects it. 
Thus, to minimize the objective \( \J_\poll^h(\cdot; 1/2) \), the components corresponding to active and idle traffic both have to be reduced, while the minimization of \( \J_\poll^h(\cdot;0) \) only requires a reduction of the contribution of active traffic.

%\begin{figure}[ht]
%	\centering
%	\input{./tikz_revision/composition_environmental_impact_review}
%	\caption{Comparison of the contribution to air pollution caused by active traffic 
%		(measured by \( \J^h_\diff \)) and idle traffic (measurement based on the average queue length).}
%	\label{fig: composition}
%\end{figure}

\Cref{fig: composition} plots the average mass of air pollutants in \( \Omega \) against the average queue length and supports these findings.
Among the three scenarios, the Pareto front to the problem~\labelcref{eq: optimization problem 3D} yields the lowest, individual contribution of active and idle traffic to air pollution. An explanation for this observation is that the optimization problem~\labelcref{eq: optimization problem 3D} tries to optimize three objectives: the traffic flow, the pollution caused by active, and the one caused by idle traffic.
Hence, contribution of active and passive traffic to air pollution has to be as low as possible instead of its sum compared to problem~\labelcref{eq: optimization problem} with \( \delta = 1/2\).
Furthermore, it becomes evident that \(\delta= 1/2 \) allows a higher traffic flow and shorter queues, while \(\delta = 0 \) allows a lower average mass of air pollutants in \( \Omega \).
The reasoning is analogous to the discussion of the results presented in \cref{fig: coordinates}.

\section{Conclusions}
We presented the framework of speed-limit dependent traffic emission models to simulate the air pollutant emission caused by vehicular and their spread given a fixed speed limit policy. Further, we used this framework together with the tool of multi-objective optimization to model our goal of achieving minimal air pollution while maximizing the traffic flow by adjusting the speed limit policy.
The proof-of-concept example presented in \cref{sec: experiments} revealed that the two objectives of minimal air pollution and maximal traffic flow are conflicting, which means they are realized for different speed limit policies. The conflict between these objectives makes multi-objective optimization a valuable tool to determine optimal compromises in the sense of Pareto where neither of the objectives is optimal.

In future work, we aim to test the framework for more complex and realistic scenarios, possibly incorporating a second-order traffic model as presented in \cite{Balzotti2022}  because these models allow to estimate the emissions based on the vehicles' speed and acceleration. 
\section*{Declarations}

\subsection*{Ethics approval and consent to participate}
Not applicable.

\subsection*{Consent for publication}
Not applicable.

\subsection*{Availability of data and materials}
Not applicable.

\subsection*{Competing interests}
The authors declare that they have no competing interests.

\subsection*{Funding}
Open Access funding enabled and organized by Projekt DEAL.

\subsection*{Authors' contributions}
S.G. and A.U. worked on the problem description and numerical simulation while M.H. provided the theoretical foundation and test instances. All authors reviewed the manuscript.

\subsection*{Acknowledgements}
S.G. thanks the Deutsche Forschungsgemeinschaft (DFG, German Research Foundation) for the financial support through GO1920/11-1 and 12-1 within the SPP 2410 Hyperbolic Balance Laws in Fluid Mechanics: Complexity, Scales, Randomness (CoScaRa).
M.H. thanks the Deutsche Forschungsgemeinschaft (DFG, German Research Foundation) for the financial support under Germany’s Excellence Strategy EXC-2023 Internet of Production 390621612 and under the Excellence Strategy of the Federal Government and the Länder, 333849990/GRK2379 (IRTG Hierarchical and Hybrid Approaches in Modern Inverse Problems), 320021702/GRK2326,
442047500/SFB1481 within the projects B04, B05 and B06, 
through SPP 2410 Hyperbolic Balance Laws in Fluid Mechanics: Complexity, Scales, Randomness (CoScaRa) within the Project(s) HE5386/26-1 and HE5386/27-1, and through SPP 2298 Theoretical Foundations of Deep Learning within the Project(s) HE5386/23-1, Meanfield Theorie zur Analysis von Deep Learning Methoden (462234017). Support through the EU DATAHYKING No. 101072546 is also acknowledged. 

\subsection*{Authors' information}
University of Mannheim, Department of Mathematics, 68131 Mannheim, Germany\\
Simone G\"ottlich and Alena Ulke\\[0.5ex]
RWTH Aachen University, Department of Mathematics, 52056 Aachen, Germany\\
Michael Herty\\[0.5ex]
{\bf Corresponding author:}\\
Correspondence to Simone G\"ottlich: {\tt goettlich@uni-mannheim.de}
\printbibliography%{}
\end{document}